\documentclass{imsart}

\RequirePackage{amsthm,amsmath,amsfonts,amssymb}
\RequirePackage[authoryear]{natbib} 
\RequirePackage[colorlinks,citecolor=blue,urlcolor=blue]{hyperref}
\RequirePackage{graphicx}

\usepackage{float,mathtools}

\startlocaldefs

\usepackage{hyperref,bbm}

\numberwithin{equation}{section}
\numberwithin{figure}{section}
\theoremstyle{plain}                    
\newtheorem{thm}{Theorem}[section]
\newtheorem{lem}{Lemma}[section]

\newtheorem{prop}{Proposition}[section]
\theoremstyle{remark}
\newtheorem{rem}{Remark}               

\def\N{{\mathbb N}}

\def\R{{\mathbb R}}

\def\1{{\mathbbm{1}}}

\def\ninfty{\mathop{\longrightarrow}\limits_{n\to\infty}}

\def\argmin{\mbox{arg}\,\min}
\newcommand{\var}{\mathop{\rm var}\nolimits}
\newcommand{\cov}{\mathop{\rm cov}\nolimits}
\newcommand{\Cov}{\mathop{\rm Cov}\nolimits}
\newcommand{\cum}{\mathop{\rm cum}\nolimits}

\newcommand{\be}{\begin{equation}}
\newcommand{\bd}{\begin{displaymath}}
\newcommand{\ed}{\end{displaymath}}
\newcommand{\bea}{\begin{eqnarray}}
\newcommand{\eea}{\end{eqnarray}}
\newcommand{\bean}{\begin{eqnarray*}}
\newcommand{\eean}{\end{eqnarray*}}

\endlocaldefs

\begin{document}

\begin{frontmatter}
\title{Trend estimation for time series with polynomial-tailed noise}
\runtitle{Trend estimation with polynomial-tailed noise}

\begin{aug}
\author[A]{\fnms{Michael H.}~\snm{Neumann}\ead[label=e1]{E-mail: michael.neumann@uni-jena.de}\orcid{0000-0002-5783-831X}}
\author[B]{\fnms{Anne}~\snm{Leucht}\ead[label=e2]{E-mail: anne.leucht@uni-bamberg.de}\orcid{0000-0003-3295-723X}}
\address[A]{Friedrich-Schiller-Universit\"at Jena,
	Institut f\"ur Mathematik,
	Ernst-Abbe-Platz 2,
	D -- 07743 Jena,
	Germany\printead[presep={,\ }]{e1}}
\address[B]{Universit\"at Bamberg,
	Institut f\"ur Statistik,
	Feldkirchenstra{\ss}e 21,
	D -- 96052 Bamberg,
	Germany\printead[presep={,\ }]{e2}}
\end{aug}

\begin{abstract}
	For time series data observed at non-random and possibly non-equidistant time points,
	we estimate the trend function nonparametrically. Under the assumption of a bounded total variation of the function 
	and low-order moment conditions on the errors we propose a nonlinear wavelet estimator which uses a Haar-type basis
	adapted to a possibly non-dyadic sample size.
	An appropriate thresholding scheme for sparse signals with an additive polynomial-tailed noise is first derived in an abstract framework
	and then applied to the problem of trend estimation.
\end{abstract}


\end{frontmatter}



\section{Introduction}
\label{S1}

\begin{figure}[h]
	\includegraphics[width=6cm]{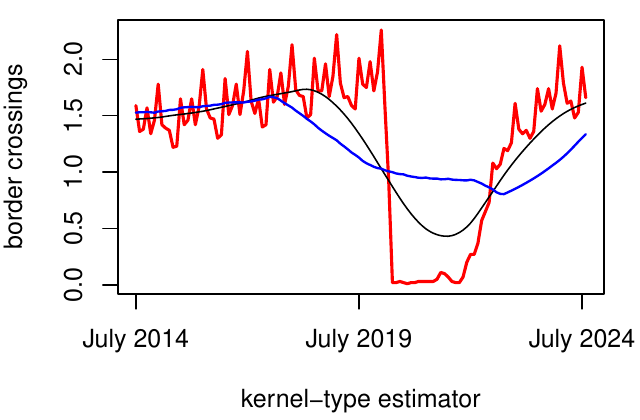}
	\caption{\textcolor{red}{\textbf{Red}}: monthly overseas arrivals in Australia (in millions) with structural breaks due to the COVID pandemic, \textbf{black}\,/\,\textcolor{blue}{\textbf{blue}}: Nadaraya Watson estimator with \textbf{Epanechnikov} and \textcolor{blue}{ \textbf{rectangular}}  kernel, bandwidth chosen by Scott's rule of thumb.}\label{fig.intro}
\end{figure}

We consider the problem of estimating the trend function of a process observed
at non-random, not necessarily equally spaced time points.
We assume that this function has a bounded total variation which includes cases
where the function is mostly smooth but has a few jumps.
In such scenarios with an inhomogeneous smoothness
a locally adapted degree of smoothing is required in order to obtain good rates
of convergence. For example, Fourier series methods, kernel estimators with a
global bandwidth and linear spline smoothers do not achieve this goal, see Figure~\ref{fig.intro} for illustration.
In contrast, a wavelet expansion of such a function provides
an efficient representation in terms of the corresponding coefficients:
only a small fraction of coefficients are large in magnitude whereas the
majority of them is small and therefore negligible. A common strategy consists 
of separating the large coefficients from the small ones by nonlinear thresholding,
which means that the former are estimated and the latter are discarded. 
Such methods have become popular in the 1990s; see for example \citet{DJKP95} and references therein.
In the case of normally distributed errors and with a sample size~$n$, a popular choice of the thresholds
is $\sqrt{2 \log(n)}$ times the standard deviation of the empirical wavelet coefficients.
\citet{Don95} showed that this choice provides a so-called ``denoising'', that is, with high probability
the estimator of the function is at least as smooth as the true function and its mean squared error
is minimax-optimal in many function classes up to a logarithmic factor.
Moreover, \citet{DJ94} proved that such estimators have a mean squared error which differs
from an ideal estimator by at most a logarithmic factor.
Even in cases with non-Gaussian and dependent errors, thresholds as in the Gaussian case can be
asymptotically appropriate since empirical wavelet coefficients at the important scales
are asymptotically normal in a sufficiently strong sense, see e.g.~\citet{Neu96}, \citet{NvS97}, and \citet{DN01}.
In the context of independent noise variables which have finite moments of a sufficiently large order,
\citet{AH05} showed that soft thresholding achieves the same asymptotic performance as in the Gaussian case.
Moreover, for independent variables with heavy tails, these authors proposed a pre-processing of the data
by median filtering which compensated for missing finite moments, and finally led to the same rate of convergence
as in the Gaussian case.
Truly non-Gaussian thresholding rules were proposed by \citet{Gao93} in the context of spectral density estimation,
by \citet{Kol99} for estimating the intensity function of a Poisson process.
In these cases, the proposed thresholds are larger than those in the Gaussian case.
Nonparametric estimation of a function with potentially inhomogeneous smoothness properties
observed at non-equidistant data points has already been considered
by a number of authors. For example, \citet{MvdG97} used least squares regression splines
regularized by a total variation penalty. 
\citet{AAdFG22} used wavelet thresholding under the assumption of sub-Gaussian noise.

In the present paper, we avoid the standard assumptions of a dyadic sample size, equally spaced sample points 
and i.i.d.~Gaussian or sub-Gaussian noise.
In Section~\ref{S2} we introduce our model and argue that the wavelet coefficients at fine scales are sparse.
More importantly, we derive our results under low-order moment conditions which are often imposed in time series analysis.
We only require finite moments up to order four
and impose a standard condition on the joint cumulants of the errors.
This results in a relatively slow polynomial decay of the tails of the distribution
of the empirical wavelet coefficients and requires an appropriate adjustment of the thresholds.
Section~\ref{S3} is devoted to a thorough discussion of estimating sparse signals
blurred by polynomial-tailed noise. It is shown in abstract Bayesian and minimax contexts
how an optimal rate of convergence depends both on the level of noise as well as the degree of sparsity
and how this rate can be attained by an appropriate choice of a threshold.
As a further prerequisite for our main results, we discuss the construction of a Haar-type basis 
for a possibly non-dyadic sample size~$n$  in Section~\ref{SS4.1}.
We use appropriately adapted basis functions which however share the essential properties
of the Haar basis. This deviates from the approach in \citet{AAdFG22} who intended to apply the standard discrete wavelet transform
and therefore proposed to embed the data points $x_1,\ldots,x_n$ into a fine equispaced grid $\{1/N,2/N,\ldots,1\}$,
where $N=2^J$ for some $J\in\N$ and $N\gg n$. 
The regularization scheme developed in Section~\ref{S3} is applied in Section~\ref{SS4.2} to our particular problem of
the estimation of a possibly discontinuous trend function of a time series.
In Section~\ref{S4+} we extend our results to a partially linear model which was also considered in \citet{AAdFG21}.
While these authors proposed to estimate the parameters of the linear part and the wavelet coefficients
simultaneously, we use a simpler approach where the linear part is first fitted by least squares
and wavelet thresholding is applied to the empirical wavelet coefficients afterwards. 
Section~\ref{S5} contains some simulations and a real data example.
It is shown that the proposed nonlinear wavelet estimator clearly outperforms kernel estimators with an optimally chosen global bandwidth.
Proofs of our main results are postponed to Section~\ref{S6} and a few auxiliary results
are collected in a final Section~\ref{S7}.

\section{Assumptions and a preview of our main results}
\label{S2}

Suppose that we observe $Y_1,\ldots,Y_n$ which form a not necessarily stationary time series and that
$EY_t=m_0(x_t)$, where $x_1<x_2<\ldots<x_n$ are ordinal variables representing e.g.~time.
This leads to the model
\begin{equation}
	\label{eq2.1}
	Y_t \,=\, m_0(x_t) \,+\, \varepsilon_t, \qquad t=1,\ldots,n.
\end{equation}
We do not assume any parametric model for the mean function~$m_0$; instead we assume that its total variation
is small in relation to the sample size~$n$.
For a real-valued function~$f$ defined on a set ${\mathcal X}\subseteq\R$, its total variation on a subset $\widetilde{\mathcal X}$ 
is defined by
${\rm TV}(f;\widetilde{\mathcal X})=\sup\big\{\sum_{i=1}^N |f(x_i)-f(x_{i-1})|\colon \;\; \{x_1,\ldots,x_N\}\subseteq\widetilde{\mathcal X},
\; x_1<x_2\ldots<x_N, \; N\in\N\big\}$.
We assume
\begin{itemize}
	\item[{\bf (A1)}\quad] $\mbox{TV}(m_0;\{x_1,\ldots,x_n\})\leq C_0$.
\end{itemize}
Regarding the errors $\varepsilon_1,\ldots,\varepsilon_n$ we impose the following weak conditions
that are standard in time series analysis.
\begin{itemize}
	\item[{\bf (A2)}\quad]
	\begin{itemize}
		\item[(i)\quad] $E\varepsilon_t=0$,
		\item[(ii)\quad] $\sup_s \sum_t \big| \cov(\varepsilon_s, \varepsilon_t) \big| \,\leq\, C_1$,
		\item[(iii)\quad] $\sup_s \sum_{t,u,v} \big| \cum(\varepsilon_s,\varepsilon_t,\varepsilon_u,\varepsilon_v) \big| \,\leq\, C_2$, 
	\end{itemize}
\end{itemize}
Here $C_0,C_1,C_2$ are fixed finite constants and\\
$\cum(\varepsilon_s,\varepsilon_t,\varepsilon_u,\varepsilon_v)\,=\,E[\varepsilon_s\,\varepsilon_t\,\varepsilon_u\,\varepsilon_v]
\,-\, E[\varepsilon_s\,\varepsilon_t]\,E[\varepsilon_u\,\varepsilon_v]
\,-\, E[\varepsilon_s\,\varepsilon_u]\,E[\varepsilon_t\,\varepsilon_v]
\,-\, E[\varepsilon_s\,\varepsilon_v]\,E[\varepsilon_t\,\varepsilon_u]$
denotes the joint cumulant of 
$\varepsilon_s,\varepsilon_t,\varepsilon_u,\varepsilon_v$.\bigskip

We propose an estimator~$\widehat{m}_n$ of $m_0$ which will be based on a wavelet expansion of this function.
Its performance is measured in terms of the mean squared error at the sample points,
i.e.~$E\big[(1/n)\sum_{t=1}^n \big(\widehat{m}_n(x_t)-m_0(x_t)\big)^2\big]$.
Under our assumption {\bf (A1)} a variant of the Haar basis is an appropriate simple choice.
We present in Section~\ref{S4} a version of this basis which is adapted to possibly unevenly spaced design points and
a non-dyadic sample size. Then the function $m_0$ can be represented as
\begin{displaymath}
	m_0(x_t) \,=\, \alpha_0^0 \quad + \quad \sum_{j=0}^{J_n}\; \sum_{k\colon\! (j,k)\in {\mathcal I}_n} \beta_{j,k}^0\, \psi_{j,k}(x_t),
\end{displaymath}
where $\alpha_0^0=(1/n)\sum_{t=1}^n m_0(x_t)$, $\beta_{j,k}^0=(1/n)\sum_{t=1}^n m_0(x_t)\,\psi_{j,k}(x_t)$,
$J_n$ such that $2^{J_n}<n\leq 2^{J_n+1}$,
and ${\mathcal I}_n\subseteq\bigcup_{j=0}^{J_n}\big(\{j\}\times\{1,\ldots,2^j\}\big)$; see Section~\ref{S4} for details.
Our wavelet estimator has the form
\begin{displaymath}
	\widehat{m}_n(x_t) \,=\, \widehat{\alpha}_0 \quad + \quad \sum_{j=0}^{J_n}\; \sum_{k\colon\! (j,k)\in {\mathcal I}_n}\widehat{\beta}_{j,k}\psi_{j,k}(x_t),
\end{displaymath}
where $\widehat{\alpha}_0$ and the $\widehat{\beta}_{j,k}$ are estimators of the corresponding coefficients.
Since the basis functions form an orthonormal system w.r.t.~the inner product $\langle f,g\rangle=(1/n)\sum_{t=1}^nf(x_t)g(x_t)$,
we have the isometry
\begin{displaymath}
	\frac{1}{n} \sum_{t=1}^n \big( \widehat{m}_n(x_t) \,-\, m_0(x_t) \big)^2
	\,=\, \sum_{j=0}^{J_n}\; \sum_{k\colon\! (j,k)\in {\mathcal I}_n} \big( \widehat{\beta}_{j,k} \,-\, \beta_{j,k}^0 \big)^2.
\end{displaymath}
This separation into the contribution of single coefficients makes an analytic study of the asymptotic performance
of $\widehat{m}_n$ possible.

Under assumption {\bf (A1)} we obtain that
\begin{equation}
	\label{eq2.2}
	2^{-j} \sum_{k\colon\! (j,k)\in {\mathcal I}_n} \big| \beta_{j,k}^0 \big| \,\leq\, C_0\, 2^{-3j/2};
\end{equation}
see (\ref{pla3.1}).
Empirical versions of the wavelet coefficients can be obtained from the sample as
$\widetilde{\alpha}_0=(1/n)\sum_{t=1}^n Y_t$ and $\widetilde{\beta}_{j,k}=(1/n)\sum_{t=1}^n m_0(x_t)Y_t$.
Under assumption {\bf (A2)} we obtain by Lemma~\ref{LA.2} that
\begin{equation}
	\label{eq2.3}
	P\big( |\widetilde{\beta}_{j,k} \,-\, \beta_{j,k}^0| \,>\, t \big)
	\,\leq\, \frac{E[ (\widetilde{\beta}_{j,k} \,-\, \beta_{j,k}^0)^4 ]}{t^4}
	\,\leq\, C\, (n^{-1/2}/t)^4 \qquad \forall t>0, \; (j,k)\in {\mathcal I}_n,
\end{equation}
for some $C<\infty$.
At fine scales~$j$, $2^{-3j/2}$ gets smaller than the noise level~$n^{-1/2}$, and the degree of sparsity
may be described by the ratio $q_{n,j}=n^{-3j/2}/n^{-1/2}$.
In the next section we study an abstract model which mimics the situation we are faced with when
we estimate the wavelet coefficients at fine scales.
This suggests how the empirical wavelet coefficients can be regularized such that the resulting estimator
attains an optimal rate of convergence.
In Section~\ref{S4} we introduce a variant of the Haar basis which is adapted to possibly unevenly spaced design points and
a non-dyadic sample size. Then we apply the regularization scheme derived in Section~\ref{S3} 
and obtain our wavelet estimator of~$m_0$.

\section{Optimal reconstruction of sparse signals from data with polynomial-tailed noise}
\label{S3}

In this section we consider an abstract model which mimics the situation we are faced with when
we estimate the wavelet coefficients at fine scales; see in particular (\ref{eq2.2}) and (\ref{eq2.3}).
Suppose first that we observe real-valued random variables~$Y_1,\ldots,Y_N$ such that
\begin{subequations}
	\begin{equation}
		\label{eq3.1a}
		Y_k \,=\, \theta_k \,+\, \varepsilon_k, \qquad k=1,\ldots,N,
	\end{equation}
	where 
	\begin{equation}
		\label{eq3.1b}
		\varepsilon_k \sim Q_\varepsilon \,\in\, {\bf Q}_\epsilon
		\,:=\, \big\{ Q\colon\quad 1-Q\big( [-t,t] \big) \,\leq\, (\epsilon/t)^4 \;\; \forall t>0 \big\}
	\end{equation}
	and 
	\begin{equation}
		\label{eq3.1c}
		\theta \,=\, \big( \theta_1,\ldots,\theta_N \big)^T \,\in\, \Theta_{N,q}
		\,:=\, \Big\{\theta\in\R^N\colon \;\;\frac{1}{N} \sum_{k=1}^N \big| \theta_k \big| \,\leq\, \epsilon\, q \Big\},
	\end{equation}
	for some $q\in (0,1)$.
\end{subequations}
Note that $\varepsilon_k$ does not have a finite fourth moment in general and it is also not required that $E\varepsilon_k=0$.
This includes the case of deterministic noise if $|\varepsilon_k|\leq\epsilon$.
In any case, it follows from Lemma~\ref{LA.1} with $t=0$ that
\begin{equation}
	\label{eq3.11}
	E\big[ \varepsilon_k^2 \big]
	\,=\, 2\, \int_0^\infty \big(1-Q_\varepsilon([-x,x])\big)\, x\, dx
	\,\leq\, 2\, \int_0^\epsilon x\, dx \;+\; 2\, \int_\epsilon^\infty \epsilon^4 x^{-3}\, dx
	\,=\, 2\, \epsilon^2,
\end{equation}
that is, $\epsilon$ may be interpreted as noise level.
The constant~$q$ in (\ref{eq3.1c}) describes the degree of sparsity of the
signals $\theta_1,\ldots,\theta_N$ in relation to~$\epsilon$.
To obtain a guideline for an appropriate regularization we consider first corresponding Bayes
and minimax problems.

To establish a Bayesian framework, we suppose that $\theta_1,\ldots,\theta_N$ are independent and follow
a three-point prior, i.e.,
\begin{subequations}
	\begin{equation}
		\label{eq3.12a}
		\theta_k \,\sim\, \pi; \qquad \pi\big( \{x\} ) \,=\, \left\{ \begin{array}{ll}
			p/2, & \quad \mbox{ if } x\in\{-\lambda, \lambda\}, \\
			1-p, & \quad \mbox{ if } x=0
		\end{array} \right. ,
	\end{equation}
	where $p=q^{4/3}$ and $\lambda=\epsilon q^{-1/3}$.
	Suppose further that the errors $\varepsilon_1,\ldots,\varepsilon_N$ are independent,
	also independent of the signals $\theta_1,\ldots,\theta_N$, and 
	\begin{equation}
		\label{eq3.12b}
		\varepsilon_k \,\sim\, \pi.
	\end{equation}
\end{subequations}
Then
\begin{displaymath}
	E\Big[ \frac{1}{N} \sum_{k=1}^N |\theta_k| \Big] \,=\, p\, \lambda \,=\, \epsilon\, q,
\end{displaymath}
that is, (\ref{eq3.1c}) is satisfied on average.
Furthermore, $P\big( |\varepsilon_k|>t \big)\leq (\epsilon/t)^4$ $\forall t>0$, that is, $\varepsilon_k \sim Q_\varepsilon\in {\bf Q}_\epsilon$.
The following results shows how the Bayes risk depends on the degree of sparsity.

{\prop
	\label{P3.1}
	Suppose that (\ref{eq3.1a}), (\ref{eq3.12a}), and (\ref{eq3.12b}) are fulfilled.
	
	Then the unique Bayes estimator $T^*(Y_k)$ is given by
	$T^*(-2\lambda)\,=\,-\lambda$, $T^*(-\lambda)\,=\,-\lambda/2$, $T^*(0)\,=\,0$, $T^*(\lambda)\,=\,\lambda/2$, $T^*(2\lambda)\,=\,\lambda$,
	and its Bayes risk is equal to
	\begin{displaymath}
		E\big[ \big( T^*(Y_k) \,-\, \theta_k \big)^2 \big] \,=\, \epsilon^2 q^{2/3}/2.
	\end{displaymath}
}
\medskip

This result can be used to obtain a lower bound to a related minimax risk.
Suppose in addition that $q^{-4/3}=O(N^{1-\gamma})$ for some $\gamma{\in(0,1)}$.
Let $\delta>0$, and let $\pi^{(N)}$ be the $N$-fold product of~$\pi$.
Let $\widetilde{T}=\widetilde{T}(Y_1,\ldots,Y_N)=(\widetilde{T}_1,\ldots,\widetilde{T}_n)^T$ be the Bayes
estimator of the vector~$\theta$ w.r.t.~the prior given by the truncation of~$\pi^{(N)}$ to $\Theta_{N,q(1+\delta)}$.
Then a lower bound to the minimax risk over $\Theta_{N,q(1+\delta)}$ is given by a corresponding Bayes risk,
\begin{eqnarray*}
	\lefteqn{ \int_{\Theta_{N,q(1+\delta)}} E_\theta \Big[ \frac{1}{N}  \sum_{k=1}^N \big( \widetilde{T}_k - \theta_k \big)^2 \Big] \, d\pi^{(N)}(\theta) } \\
	& = & \int_{\R^N} E_\theta \Big[ \frac{1}{N}  \sum_{k=1}^N \big( \widetilde{T}_k - \theta_k \big)^2 \Big] \, d\pi^{(N)}(\theta)
	\,-\, \int_{\R^N\setminus \Theta_{N,q(1+\delta)}} E_\theta \Big[\frac{1}{N}  \sum_{k=1}^N \big( \widetilde{T}_k - \theta_k \big)^2 \Big] \, d\pi^{(N)}(\theta).
\end{eqnarray*}
The first term on the right-hand side can be estimated from below by\\
$\int_{\R^N} E_\theta \big[ (1/N) \sum_{k=1}^N \big( T^*(Y_k) - \theta_k \big)^2 \big] \, d\pi^{(N)}(\theta) \,=\, \epsilon^2 q^{2/3}/2$.
Regarding the second one, note first that $\widetilde{T}_k\in [-\lambda,\lambda]$ holds with probability~1,
which leads to
\begin{displaymath}
	E_\theta \Big[\frac{1}{N}  \sum_{k=1}^N \big( \widetilde{T}_k - \theta_k \big)^2 \Big] \,\leq\, (2\lambda)^2 \,=\, 4\epsilon^2 q^{-2/3}.
\end{displaymath}
On the other hand, we obtain from Bernstein's inequality that
\begin{eqnarray*}
	P\big( \theta \not\in \Theta_{N,q(1+\delta)} \big) 
	& = & P\Big( \frac{1}{N}  \sum_{k=1}^N |\theta_k| - E|\theta_k| \,>\, \epsilon\, q (1+\delta) \Big) \\
	& \leq & \exp\bigg\{ - \frac{ (N\epsilon q(1+\delta))^2/2 }
	{ \sum_{k=1}^N E[ (|\theta_k|-E|\theta_k|)^2 ] \;+\; \big\| |\theta_1|-E|\theta|_1 \big\|_\infty \, (N\epsilon q(1+\delta))/3 } \bigg\} \\
	& = & \exp\big\{ - R_N \big\},
\end{eqnarray*}
where 
\begin{displaymath}
	1/ R_N \,\asymp\, \frac{ N\lambda^2p \;+\; \lambda N\epsilon q }{ (N\epsilon q)^2 } \,=\, O\big( N^{-1} q^{-4/3} \big)
	\,=\, O\big( N^{-\gamma} \big).
\end{displaymath}
This implies that the second term is of order $O\big( \epsilon^2 q^{-2/3} \exp\{-R_N\} \big)$ and we obtain that
\begin{equation}
	\label{eq3.3}
	\begin{aligned}
	&\inf_{\widehat{T}} \sup\bigg\{ E_\theta\Big[ \frac{1}{N} \sum_{k=1}^N (\widehat{T}_k-\theta)^2 \Big]\colon \quad
	\theta\in \Theta_{N,q(1+\delta)}, \; \varepsilon_1,\ldots,\varepsilon_N\sim Q_\varepsilon\in {\bf Q}_\epsilon \bigg\}\\
	&\,\geq\, \big(\epsilon^2 \, q^{2/3} /2\big) (1+o(1)),
	\end{aligned}
\end{equation}
as $N\to\infty$.

Although the form of~$T^*$ might suggest that the linear estimator~$Y_k/2$ of~$\theta_k$
is appropriate, this is no longer true for some other error distributions which satisfy (\ref{eq3.1b}).
If, for example, $\varepsilon_k\sim {\mathcal N}(0,\epsilon^2)$, then it follows from
$E\big[ (Y_k/2-\theta_k)^2 \big]\,=\,\epsilon^2/4+(\theta_k/2)^2$ that the Bayes risk of~$Y_k/2$ under~$\pi$
is equal to $\epsilon^2/4+p\lambda^2/4=\epsilon^2(1+q^{2/3})/4$, that is,
it is dominated by the variance of $Y_k$ and there is not the desired gain due to the sparsity of the signal~$\theta_k$.
This is of course in line with common folklore that nonlinear methods are required in order to 
efficiently estimate sparse signals.

We therefore consider estimators $\widehat{\theta}_k=\widehat{T}(Y_k)$ which satisfy the following conditions:
\begin{subequations}
	\begin{equation}
		\label{eq3.3a}
		\widehat{T}(y) \,=\, 0, \qquad \mbox{ if } |y|\,<\, t
	\end{equation}
	and
	\begin{equation}
		\label{eq3.3b}
		\big| \widehat{T}(y) \,-\, y \big| \,\leq\, t, \qquad \mbox{ if } |y|\,\geq\, t.
	\end{equation}
\end{subequations}
Note that the popular methods of soft and hard thresholding with threshold~$t$ satisfy
these conditions.
The following theorem shows that the above estimator attains the optimal rate of convergence.

{\thm
	\label{T3.1}
	Suppose that (\ref{eq3.1a}), (\ref{eq3.1b}), (\ref{eq3.3a}), and (\ref{eq3.3b}) are fulfilled.
	Then
	\begin{itemize}
		\item[(i)\quad] $E\big[ (\widehat{\theta}_k - \theta_k)^2 \big]
		\,\leq\, \left\{ \begin{array}{ll}
			\theta_k^2 \,+\, 9E[ \varepsilon_k^2 \, \1(|\varepsilon_k|>t/2) ] & \qquad \mbox{ if } |\theta_k|\leq t/2, \\
			2t^2 \,+\, 4\epsilon^2 & \qquad \mbox{ if } |\theta_k|> t/2
		\end{array} \right.$ .
		\item[(ii)\quad] If in addition $t=K\lambda=K\epsilon q^{-1/3}$
		for some $K>0$, then
		\begin{displaymath}
			\sup_{\theta\in\Theta_{N,q(1+\delta)},\, \varepsilon_1,\ldots,\varepsilon_N\sim Q_\varepsilon\in {\bf Q}_\epsilon}
			E_\theta \Big[  \frac{1}{N}\sum_{k=1}^N \big( \widehat{\theta}_k \,-\, \theta_k \big)^2 \Big]
			\,=\, O\big( \epsilon^2\, q^{2/3} \big).
		\end{displaymath}
	\end{itemize}
}

\section{Main results}
\label{S4}

\subsection{A Haar-type basis for unevenly spaced data}
\label{SS4.1}

Let us first consider the simple case where we have to deal with a real-valued function~$f$ on the 
unit interval~(0,1]. This interval can be decomposed into intervals
\begin{displaymath}
	I_{j,k} \,=\, \big( (k-1)2^{-j}, k2^{-j} \big] ,\qquad  j=0,1,2,\ldots;\; k=1,\ldots,2^j.
\end{displaymath}
To define the Haar basis we start with a scaling function $\phi_0 \,=\, \1_{I_{0,1}} \,=\, \1_{(0,1]}$
and a mother wavelet $\psi \,=\, \1_{I_{1,1}}-\1_{I_{1,2}} \,=\, \1_{(0,1/2]}-\1_{(1/2,1]}$.
Using dilations and translations we obtain wavelets
$\psi_{j,k} \,=\, 2^{j/2}\,\big( \1_{I_{j+1,2k-1}}-\1_{I_{j+1,2k}} \big)
\,=\, 2^{j/2}\, \psi(2^j\cdot \,-\, (k-1))$, where $j\geq 0;\; k=1,\ldots,2^j$.
The collection of these functions $\big\{\phi_0, \psi_{j,k} (j\geq 0; k=1,\ldots,2^j)\big\}$ forms an
orthonormal basis of $L_2((0,1])$. An arbitrary function $f\in L_2((0,1])$ can be expanded as
\begin{displaymath}
	f(x) \,=\, \alpha_0\, \phi_0(x) \,+\, \sum_{j=0}^\infty \sum_{k=1}^{2^j} \beta_{j,k}\, \psi_{j,k}(x),
\end{displaymath}
where $\alpha_0\,=\,\int \phi_0(x) f(x)\, dx$ and
$\beta_{j,k}\,=\,\int \psi_{j,k}(x) f(x)\, dx$.
Since the wavelets $\psi_{j,1},\ldots,\psi_{j,2^j}$ are supported on respective disjoint intervals $I_{j,1},\ldots,I_{j,2^j}$,
we obtain the following estimate for the size of the wavelet coefficients:
\begin{equation}
	\label{eq4.0}
	\sum_{k=1}^{2^j} \big|\beta_{j,k}\big| \,=\, \sum_{k=1}^{2^j} \|\psi_{j,k}\|_1 \, \cdot \, \big( {\rm TV}(f;I_{j,k})/2 \big)
	\,\leq\, 2^{-j/2-1}\,{\rm TV}(f;(0,1]).
\end{equation}

Suppose for the time being that we have a dyadic sample size $n=2^{J_n+1}$ for some $J_n\in\N$
and that $x_t=t/n$ for all $t=1,\ldots,n$.
Then we can directly use the Haar basis for our purposes.
Let $\bar{m}_0$ be a piecewise constant continuation of $m_0(x_1),\ldots,m_0(x_n)$ on the intervals $(x_{t-1},x_t]$,
i.e., $\bar{m}_0(x)=m_0(x_t)$ $\forall x\in (x_{t-1},x_t]$.
Then $\bar{m}_0$ can be perfectly represented by the first~$n$ basis functions,
\begin{displaymath}
	\bar{m}_0(x) \,=\, \alpha_0^0\, \phi_0(x) \,+\, \sum_{j=0}^{J_n} \sum_{k=1}^{2^j} \beta_{j,k}^0\, \psi_{j,k}(x),
\end{displaymath}
where $\alpha_0^0 \,=\, \int \phi_0(x)\, \bar{m}_0(x) \, dx \,=\, (1/n) \sum_{t=1}^n m_0(x_t)$
and $\beta_{j,k}^0 \,=\, \int \psi_{j,k}(x)\, \bar{m}_0(x) \, dx \,=\, (1/n) \sum_{t=1}^n \psi_{j,k}(x_t)\, m_0(x_t)$.
We have in particular that
\begin{displaymath}
	m_0(x_t) \,=\, \alpha_0^0 \,+\, \sum_{j=0}^{J_n}\, \sum_{k=1}^{2^j} \beta_{j,k}^0\, \psi_{j,k}(x_t) \qquad \forall t=1,\ldots,n.
\end{displaymath}
The observations $Y_1,\ldots,Y_n$ can also be described this way,
\begin{displaymath}
	Y_t \,=\, \widetilde{\alpha}_0\phi_0(x_t) \,+\, \sum_{j=0}^{J_n} \sum_{k=1}^{2^j} \widetilde{\beta}_{j,k}\, \psi_{j,k}(x_t),
\end{displaymath}
where $\widetilde{\alpha}_0\,=\,(1/n) \sum_{t=1}^n \phi_0(x_t) Y_t\,=\,\bar{Y}_n$ and
$\widetilde{\beta}_{j,k}\,=\,(1/n) \sum_{t=1}^n \psi_{j,k}(x_t) Y_t$ are empirical versions of the respective wavelet coefficients.
Since the vectors $\big(\phi_0(x_1),\ldots,\phi_0(x_n)\big)^T$ and $\big(\psi_{j,k}(x_1),\ldots,\psi_{j,k}(x_n)\big)^T$ 
($j=0,\ldots,J_n; \; k=1,\ldots,2^j$) form an orthonormal basis of~$\R^n$ w.r.t.~the inner product $\langle\cdot,\cdot\rangle$
given by $\langle a,b\rangle=(1/n)\sum_{t=1}^n a_i b_i$, we also obtain,
for $\widehat{m}_n(x)=\widehat{\alpha}_0\,+\,\sum_{j=0}^{J_n} \sum_{k=1}^{2^j} \widehat{\beta}_{j,k} \, \psi_{j,k}(x)$, that
\begin{displaymath}
	\frac{1}{n} \sum_{t=1}^n \big( \widehat{m}_n(x) \,-\, m_0(x_t) \big)^2
	\,=\, \big( \widehat{\alpha}_0 \,-\, \alpha_0^0 \big)^2 
	\,+\, \sum_{j=0}^{J_n} \sum_{k=1}^{2^j} \big( \widehat{\beta}_{j,k} \,-\, \beta_{j,k}^0 \big)^2.
\end{displaymath}

If the sample size~$n$ is not a dyadic one, then we have to modify the above approach slightly.
Since the~$x_t$ can be considered as ordinal variables, we simplify our notation by sticking to our assumption that $x_t=t/n$.
We choose the finest scale~$J_n$ such that
\begin{displaymath}
	2^{J_n} \,<\, n \,\leq\, 2^{J_n+1}.
\end{displaymath}
As above, we take as a starting point dyadic intervals on $(0,1]$:
\begin{displaymath}
	I_{j,k} \,=\, \big( (k-1)2^{-j}, k2^{-j} \big], \qquad j=0,\ldots,J_n+1;\, k=1,\ldots,2^j.
\end{displaymath}
It suffices to the specify the functions $\phi_0$ and $\psi_{j,k}$ at the points $x_1,\ldots,x_n$.
Let
\begin{displaymath}
	\phi_0(x_t) \,=\, 1 \qquad \forall t=1,\ldots,n.
\end{displaymath}
We define wavelet functions $\psi_{j,k}$ such that the essential properties of the Haar basis are retained:
These functions shall form an orthonormal system w.r.t.~the inner product $\langle\cdot,\cdot\rangle$
and they should efficiently describe functions with jumps.
Let
\begin{displaymath}
	n_{j,k} \,:=\, \#\big\{ 1\leq t\leq n\colon\; x_t\in I_{j,k} \big\}.
\end{displaymath}
It is not difficult to see that
\begin{displaymath}
	[n 2^{-j}] \,\leq\, n_{j,k} \,<\, n2^{-j} \,+\, 1.
\end{displaymath}
(Since the length of the interval $I_{j,k}$ is $2^{-j}$ and the distance between adjacent points $x_{t-1}$ and $x_t$
is $1/n$, we obtain that $n_{j,k}=n2^{-j}$ if $n2^{-j}$ is an integer.
Otherwise, if $n2^{-j}$ is not an integer, then $n_{j,k}\geq [n2^{-j}]$.
On the other hand, $n_{j,k}\geq n2^{-j}+1$ is impossible since this implies $n_{j,k}\geq [n2^{-j}]+2$ and so
the length of $I_{j,k}$ would exceed $([n2^{-j}]+1)/n$. This, however, leads to a contradiction since
$([n2^{-j}]+1)/n>n2^{-j}/n=2^{-j}$.)
First we obtain a system of orthogonal functions by
\begin{displaymath}
	\widetilde{\psi}_{j,k} \,=\, \frac{1}{n_{j+1,2k-1}} \1_{I_{j+1,2k-1}} \,-\, \frac{1}{n_{j+1,2k}} \1_{I_{j+1,2k}}
	\qquad \forall (j,k)\in {\mathcal I}_n,
\end{displaymath} 
where ${\mathcal I}_n:=\big\{ (j,k)\colon\; n_{j+1,2k-1}\geq 1 \mbox{ and } n_{j+1,2k}\geq 1 \big\}$.
Since $n_{J_n+1,k}<n2^{-(J_n+1)}+1\leq 2$, we obtain $n_{J_n+1,k}\leq 1$ for all $k=1,\ldots,2^{J_n+1}$.
This implies that $(j,k)\not\in {\mathcal I}_n$, if $j>J_n$, i.e., $J_n$ is the finest scale where functions $\widetilde{\psi}_{j,k}$
are defined.
Moreover, the `decomposition pyramid' does not stop before the interval $(0,1]$ is decomposed into the smallest possible
intervals $((t-1)/n,t/n]$. For example, if $n_{j,k}=2$, then $n_{j+1,2k-1}=n_{j+1,2k}=1$, and so $(j,k)\in {\mathcal I}_n$.
This implies in particular that $\# {\mathcal I}_n=n-1$.

We have
\begin{displaymath}
	\langle \phi_0, \widetilde{\psi}_{j,k} \rangle \,:=\, \frac{1}{n} \sum_{t=1}^n \phi_0(x_t) \, \widetilde{\psi}_{j,k}(x_t) \,=\, 0
	\qquad \forall (j,k)\in {\mathcal I}_n
\end{displaymath}
and
\begin{displaymath}
	\langle \widetilde{\psi}_{j,k}, \widetilde{\psi}_{j',k'} \rangle
	\,:=\, \frac{1}{n} \sum_{t=1}^n \widetilde{\psi}_{j,k}(x_t) \, \widetilde{\psi}_{j',k'}(x_t) \,=\, 0
	\qquad \forall (j,k),(j',k')\in {\mathcal I}_n,\, (j,k)\neq (j',k'),
\end{displaymath}
i.e., these functions form an orthogonal system w.r.t.~the inner product $\langle\cdot,\cdot\rangle$.
It remains to normalize these functions.
Since
\begin{displaymath}
	\langle \widetilde{\psi}_{j,k}, \widetilde{\psi}_{j,k} \rangle
	\,=\, \frac{1}{n} \big( \frac{1}{n_{j+1,2k-1}} \,+\, \frac{1}{n_{j+1,2k}} \big), 
\end{displaymath}
we obtain by
\begin{displaymath}
	\psi_{j,k} \,:=\, \widetilde{\psi}_{j,k}/\sqrt{ \langle \widetilde{\psi}_{j,k}, \widetilde{\psi}_{j,k} \rangle }
	\,=\, \frac{ \sqrt{n} }{ \sqrt{ \frac{1}{n_{j+1,2k-1}} \,+\, \frac{1}{n_{j+1,2k}} } }
	\; \big( \frac{1}{n_{j+1,2k-1}} \1_{I_{j+1,2k-1}} \,-\, \frac{1}{n_{j+1,2k}} \1_{I_{j+1,2k}} \big)
\end{displaymath}
that $\langle \psi_{j,k}, \psi_{j,k}\rangle=1$.
The vectors
\begin{displaymath}
	\boldsymbol{\phi}_0 \,=\, (1/\sqrt{n}) \big( \phi_0(x_1),\ldots,\phi_0(x_n) \big)^T
\end{displaymath}
and
\begin{displaymath}
	\boldsymbol{\psi}_{j,k} \,=\, (1/\sqrt{n}) \big( \psi_{j,k}(x_1),\ldots,\psi_{j,k}(x_n) \big)^T \qquad \forall (j,k)\in {\mathcal I}_n
\end{displaymath}
form an othonormal system in $\R^n$. Since $\#{\mathcal I}_n=n-1$, this is even an orthonormal basis.
As in the case of $n=2^{J_n+1}$ we have,
for $\widehat{m}_n(x)=\widehat{\alpha}_0\,+\,\sum_{j=0}^{J_n} \; \sum_{k\colon\! (j,k)\in{\mathcal I}_n} \widehat{\beta}_{j,k} \, \psi_{j,k}(x)$,
the isometry
\begin{equation}
	\label{eq41.2}
	\frac{1}{n} \sum_{t=1}^n \big( \widehat{m}_n(x) \,-\, m_0(x_t) \big)^2
	\,=\, \big( \widehat{\alpha}_0 \,-\, \alpha_0^0 \big)^2 
	\,+\, \sum_{j=0}^{J_n} \; \sum_{k\colon\! (j,k)\in{\mathcal I}_n} \big( \widehat{\beta}_{j,k} \,-\, \beta_{j,k}^0 \big)^2.
\end{equation}

\subsection{A nonlinear wavelet estimator of the trend function}
\label{SS4.2} 

Now we consider an estimator $\widehat{m}_n$ of $m_0$, where
\begin{displaymath}
	\widehat{m}_n(x_t) \,=\, \widehat{\alpha}_0 \,+\, \sum_{j=0}^{J_n} \sum_{k\colon\! (j,k)\in {\mathcal I}_n} \widehat{\beta}_{j,k}\, \psi_{j,k}(x_t).
\end{displaymath}
The coefficients $\widehat{\alpha}_0$ and $\widehat{\beta}_{j,k}$ of this wavelet expansion are derived from corresponding empirical versions
$\widetilde{\alpha}_0=(1/n)\sum_{t=1}^n \phi_0(x_t) Y_t$
and $\widetilde{\beta}_{j,k}=(1/n)\sum_{t=1}^n \psi_{j,k}(x_t) Y_t$ of the true coefficients
$\alpha_0^0\,=\,(1/n) \sum_{t=1}^n \phi_0(x_t) m_0(x_t)$ and
$\beta_{j,k}^0\,=\,(1/n) \sum_{t=1}^n \psi_{j,k}(x_t) m_0(x_t)$ of the function~$m_0$.
In view of the isometry (\ref{eq41.2}) above,
we direct our attention to the estimation of the coefficients.
It follows from Lemma~\ref{LA.2} that
\begin{equation}
	\label{eq4.11}
	E\big[ \big( \widetilde{\beta}_{j,k} \,-\, \beta_{j,k}^0 \big)^4 \big] \,\leq\, C_3\, n^{-2} \qquad \forall (j,k)\in {\mathcal I}_n,
\end{equation}
for some $C_3<\infty$, which implies that
\begin{equation}
	\label{eq4.12}
	P\big( |\widetilde{\beta}_{j,k} \,-\, \beta_{j,k}^0|>t \big) \,\leq\, C_3\, (n^{-1/2}/t)^4 \qquad \forall t>0.
\end{equation}
On the other hand, we obtain similarly to (\ref{eq4.0}) that
\begin{equation}
	\label{eq4.13}
	2^{-j} \sum_{k\colon\! (j,k)\in {\mathcal I}_n} |\beta_{j,k}^0| \,\leq\, C_0\, 2^{-3j/2}.
\end{equation}
This means that the signal becomes sparse in relation to the noise level~$n^{-1/2}$ at scales~$j$
where $2^{-3j/2}$ is of smaller order than $n^{-1/2}$.
The degree of sparsity is expressed by $q_{n,j}=n^{1/2}2^{-3j/2}$.
In view of the message provided by Theorem~\ref{T3.1} we will modify the empirical coefficients $\widetilde{\beta}_{j,k}$
nonlinearly at scales $j\geq J_n^*$, where $J_n^*$ is the critical level.
Let, for definiteness, $J_n^*$ be such that $2^{J_n^*-1}<n^{1/3}\leq 2^{J_n^*}$ and let $t_{n,j}=K\, n^{-1/2} q_{n,j}^{-1/3}=K\,n^{-2/3}2^{j/2}$,
where $K$ is an arbitrary positive constant.
We focus on the estimator
\begin{displaymath}
	\widehat{m}_n(x_t) \,=\, \widetilde{\alpha}_0
	\,+\, \sum_{j=0}^{J_n^*-1} \; \sum_{k\colon\! (j,k)\in {\mathcal I}_n} \widetilde{\beta}_{j,k} \, \psi_{j,k}(x_t)
	\,+\, \sum_{j=J_n^*}^{J_n} \; \sum_{k\colon\! (j,k)\in {\mathcal I}_n} \widehat{\beta}_{j,k} \, \psi_{j,k}(x_t),
\end{displaymath}
where, for $(j,k)\in {\mathcal I}_n$, $j\geq J_n^*$,
\begin{subequations}
	\begin{equation}
		\label{eq4.1a}
		\widehat{\beta}_{j,k} \,=\, 0, \qquad \mbox{ if } |\widetilde{\beta}_{j,k}|< t_{n,j}
	\end{equation}
	and
	\begin{equation}
		\label{eq4.1b}
		\big|\widehat{\beta}_{j,k} \,-\, \widetilde{\beta}_{j,k}\big| \,\leq\, t_{n,j}, \qquad \mbox{ if } |\widetilde{\beta}_{j,k}|\geq t_{n,j}.
	\end{equation}
\end{subequations}
Popular examples of such a strategy are hard thresholding,
\begin{displaymath}
	\widehat{\beta}_{j,k}^{(h)} \,=\, \left\{ \begin{array}{ll} 0, & \qquad \mbox{ if } |\widetilde{\beta}_{j,k}| \,<\, t_{n,j}, \\
		\widetilde{\beta}_{j,k}, & \qquad \mbox{ if } |\widetilde{\beta}_{j,k}| \,\geq\, t_{n,j}
	\end{array} \right. ,
\end{displaymath}
and soft thresholding,
\begin{displaymath}
	\widehat{\beta}_{j,k}^{(s)} \,=\, \mbox{sgn}(\widetilde{\beta}_{j,k})\, \big( |\widetilde{\beta}_{j,k}| \,-\, t_{n,j} \big)^+.
\end{displaymath}
Note, in passing, that there is a well-known connection between soft thresholding and $l_1$-penalization,
i.e., $\widehat{\beta}_{j,k}^{(s)}$ is the unique minimizer of the function $\beta\mapsto (\beta-\widetilde{\beta}_{j,k})^2+2t_{n,j}|\beta|$.
It follows from (\ref{eq4.11}) that
\begin{equation}
	\label{eq4.2a}
	E\big[ \big( \widetilde{\alpha}_0 \,-\, \alpha_0^0 \big)^2 \big]
	\,+\, \sum_{j=0}^{J_n^*-1}\; \sum_{k\colon\! (j,k)\in {\mathcal I}_n} E\big[ \big( \widetilde{\beta}_{j,k} \,-\, \beta_{j,k}^0 \big)^2 \big]
	\,=\, O\big( 2^{J_n^*}\, n^{-1} \big) \,=\, O\big( n^{-2/3} \big)
\end{equation}
In order to estimate the contribution to the risk by the other coefficients we use the following result.
It improves the simple upper estimate $\sum_{k\colon\! (j,k)\in {\mathcal I}_n} (\beta_{j,k}^0)^2\wedge t_{n,j}^2
\leq t_{n,j}\; \sum_{k\colon (j,k)\in {\mathcal I}_n} |\beta_{j,k}^0|=O(n^{-2/3})$ and it is essential for the proof of Theorem~\ref{T4.2} below.

{\lem
	\label{L4.1}
	Let $m_0\colon (0,1]\rightarrow\R$ be such that
	$\mbox{TV}\big( m_0; \{x_1,\ldots,x_n\} \big) \,\leq\, C_0$ and let $t_{n,j}=Kn^{-2/3}2^{j/2}$.
	Then, for $\beta_{j,k}^0=(1/n)\sum_{t=1}^n \psi_{j,k}(x_t) m(x_t)$,
	\begin{displaymath}
		S_n(m_0) \,:=\, \sum_{j\colon 2^j\geq n^{1/3}}\; \sum_{k\colon\! (j,k)\in {\mathcal I}_n} (\beta_{j,k}^0)^2 \wedge t_{n,j}^2
		\,=\, O\big( n^{-2/3} \big).
	\end{displaymath}
}
\medskip

Since $\var(\widetilde{\beta}_{j,k})=O(t_{n,j})$ for all $(j,k)\in{\mathcal I}_n$, $j\geq J_n^*$, we obtain
from (i) of Theorem~\ref{T3.1}, that for $j\geq J_n^*$,
\begin{eqnarray*}
	\lefteqn{ \sum_{k\colon\! (j,k)\in {\mathcal I}_n} E\big[ (\widehat{\beta}_{j,k} - \beta_{j,k}^0)^2 \big] } \\
	& = & \sum_{k\colon\! (j,k)\in {\mathcal I}_n} O\Big( \min\big\{ (\beta_{j,k}^0)^2, t_{n,j}^2 \big\}
	\;+\; \, E\big[ (\widetilde{\beta}_{j,k}-\beta_{j,k}^0)^2 \, \1\big( |\widetilde{\beta}_{j,k}-\beta_{j,k}^0|>t_{n,j}/2 \big) \big] \Big).
\end{eqnarray*}
Since
\begin{displaymath}
	\sum_{k\colon\! (j,k)\in {\mathcal I}_n} \min\big\{ (\beta_{j,k}^0)^2, t_{n,j}^2 \big\}
	\,\leq\, t_{n,j} \; \sum_{k\colon\! (j,k)\in {\mathcal I}_n} |\beta_{j,k}^0| \,=\, O\big( n^{-2/3} \big)
\end{displaymath}
and 
\begin{displaymath}
	E\big[ (\widetilde{\beta}_{j,k}-\beta_{j,k}^0)^2 \, \1\big( |\widetilde{\beta}_{j,k}-\beta_{j,k}^0|>t_{n,j}/2 \big) \big]
	\,\leq\, E\big[ (\widetilde{\beta}_{j,k}-\beta_{j,k}^0)^4 \big]/(t_{n,j}/2)^2 \,=\,  O\big( n^{-2/3} 2^{-j} \big),
\end{displaymath}
it follows that
\begin{equation}
	\label{eq4.2b}
	\sum_{k\colon\! (j,k)\in {\mathcal I}_n} E\big[ \big( \widehat{\beta}_{j,k} \,-\, \beta_{j,k}^0 \big)^2 \big]
	\,=\, O\big( n^{-2/3} \big).
\end{equation}
Since $J_n-J_n^*=O(\log(n))$, we obtain from (\ref{eq4.2a}) and (\ref{eq4.2b}) the following result.
\bigskip

{\thm
	\label{T4.1}
	Suppose that {\bf (A2)} is fulfilled. Then
	\begin{displaymath}
		\sup\Big\{ E\Big[ \frac{1}{n} \sum_{t=1}^n \big( \widehat{m}_n(x_t) \,-\, m_0(x_t) \big)^2 \Big]\colon \;\; TV(m_0;[0,1])\leq C_0 \Big\}
		\,=\, O\big( n^{-2/3}\, \log(n) \big).
	\end{displaymath}
}

This result can be improved if we replace condition {\bf (A2)} by the following majorization condition for the distribution
of the empirical wavelet coefficients.

\begin{itemize}
	\item[{\bf (A2')}\quad] 
	\begin{itemize}
		\item[(i)\quad]
		$E\big[ \big( \widetilde{\alpha}_0 \,-\, \alpha_0^0 \big)^2 \big]
		\,+\, \sum_{j=0}^{J_n^*-1} \; \sum_{k\colon\! (j,k)\in {\mathcal I}_n} E\big[ \big( \widetilde{\beta}_{j,k} \,-\, \beta_{j,k}^0 \big)^2 \big]
		\,=\, O\big( n^{-2/3} \big).$
		\item[(ii)\quad]
		There exists a distribution function~$G$ on $[0,\infty)$ (not necessarily a probability
		distribution function) such that $\int_0^\infty x^4 \, dG(x) <\infty$ and
		\begin{equation}
			\label{eq4.3}
			P\big( n^{1/2} \big| \widetilde{\beta}_{j,k}-\beta_{j,k}^0 \big| > x \big) \,\leq\, 1 \,-\, G(x)
			\qquad \forall x\geq 0
		\end{equation}
		holds for all $(j,k)\in {\mathcal I}_n$, $j\geq J_n^*$.
	\end{itemize}
\end{itemize}

{\rem
	Condition (\ref{eq4.3}), as it stands, is a high-level condition and such conditions should be avoided wherever possible.
	There are however scenarios where this relation follows from simple conditions on the errors $\varepsilon_1,\ldots,\varepsilon_n$.
	\begin{itemize}
		\item[1)\quad]
		If $\varepsilon_1,\ldots,\varepsilon_n$ are jointly normal with zero mean, then
		\begin{displaymath}
			n^{1/2} \big( \widetilde{\beta}_{j,k}-\beta_{j,k}^0 \big) \,\sim\, N\big( 0, \sigma_{j,k}^2 \big),
		\end{displaymath}
		where
		\begin{eqnarray*}
			\sigma_{j,k}^2 & = & n^{-1} \sum_{s,t=1}^n \psi_{j,k}(x_s) \psi_{j,k}(x_t) \cov( \varepsilon_s, \varepsilon_t )
			\,\leq\, n^{-1} \sum_{s=1}^n \psi_{j,k}^2(x_s) \; \sum_{t=1}^n \big| \cov( \varepsilon_s, \varepsilon_t ) \big| \\
			& \leq & \bar{\sigma}^2 \,:=\, \max_{1\leq s\leq n} \sum_{t=1}^n \big| \cov( \varepsilon_s, \varepsilon_t ) \big|.
		\end{eqnarray*}
		In this case, $G(x)=2\Phi(x/\bar{\sigma})-1$ satisfies (\ref{eq4.3}).
		\item[2)\quad]
		If $\varepsilon_1,\ldots,\varepsilon_n$ are independent with zero mean
		\begin{displaymath}
			\max_{1\leq t\leq n} E\big[ |\varepsilon_t|^\gamma \big] \,<\, \infty
		\end{displaymath}
		for some $\gamma>4$, 
		then we obtain by Rosenthal's inequality (see Theorem~3 in \citet{Ros70}) that
		\begin{eqnarray*}
			\lefteqn{ E\big[ \big| \sqrt{n}(\widetilde{\beta}_{j,k}-\beta_{j,k}^0) \big|^\gamma \big] } \\
			& = & E\big[ \big| \frac{1}{\sqrt{n}} \sum_{t=1}^n \psi_{j,k}(x_t) \varepsilon_t \big|^\gamma \big] \\
			& \leq & C_\gamma \, \max\Big\{ \sum_{t=1}^n E\big[ \big| \psi_{j,k}(x_t)\varepsilon_t/\sqrt{n} \big|^\gamma \big],
			\Big( \sum_{t=1}^n E\big[ \big| \psi_{j,k}(x_t)\varepsilon_t/\sqrt{n} \big|^2 \Big)^{\gamma/2} \Big\} \\
			& \leq & C_\gamma \, \max\big\{ E\big[ |\varepsilon_t|^\gamma \big]\colon \; 1\leq t\leq n \big\} \,=:\, C_\gamma'.
		\end{eqnarray*}
		(The latter inequality follows from $(\psi_{j,k}(x_t)/\sqrt{n})^2\leq (1/n)\sum_{t=1}^n \psi_{j,k}^2(x_t)=1$.)
		Then
		\begin{displaymath}
			P\big( n^{1/2} \big| \widetilde{\beta}_{j,k}-\beta_{j,k}^0 \big| > x \big)
			\,\leq\, \min\big\{ 1, C_\gamma'/x^\gamma \big\} \,=:\, 1 \,-\, G(x) \qquad \forall x\geq 0.
		\end{displaymath}
		and it holds that $\int_0^\infty x^4\, dG(x)<\infty$, as required.
		\item[3)\quad]
		If, for $\gamma>4$, $\epsilon>0$, and some even number $c>\gamma$,
		$\varepsilon_1,\ldots,\varepsilon_n$ are strong ($\alpha-$) mixing with coefficients satisfying
		\begin{displaymath}
			\sum_{r=1}^\infty (r+1)^{c-2} \big( \alpha(r) \big)^{\epsilon/(c+\epsilon)} \,<\, \infty,                                
		\end{displaymath}
		and if the $\varepsilon_t$ have zero mean and $\max_{1\leq t\leq n} E\big[ |\varepsilon_t|^{\gamma+\epsilon} \big]<\infty$,
		then we obtain from a Rosenthal-type inequality (see e.g.~\citet[page~26, Theorem~2]{Dou94}) that
		\begin{displaymath}
			E\big[ \big| \sqrt{n}(\widetilde{\beta}_{j,k}-\beta_{j,k}^0) \big|^\gamma \big]
			\,=\, E\big[ \big| \frac{1}{\sqrt{n}} \sum_{t=1}^n \psi_{j,k}(x_t) \varepsilon_t \big|^\gamma \big]
			\,\leq\, C\, \Big( E\big[ |\varepsilon_t|^{\gamma+\epsilon} \big] \Big)^{\gamma/(\gamma+\epsilon)}
			\,=:\, C_{\gamma}^{''}.
		\end{displaymath}
		Then~$G$ given by $G(x):= 1-\min\big\{ 1, C_\gamma^{''}/x^\gamma \big\}$ satisfies (\ref{eq4.3}).
	\end{itemize}
}
\bigskip

\noindent
Under {\bf (A2')} we obtain a rate of convergence without a logarithmic factor which is known to be optimal
in many similar estimation problems.

{\thm
	\label{T4.2}
	If {\bf (A2')} is fulfilled, then
	\begin{displaymath}
		\sup\Big\{ E\Big[ \frac{1}{n} \sum_{t=1}^n \big( \widehat{m}_n(x_t) \,-\, m_0(x_t) \big)^2 \Big]\colon \;\;
		TV(m_0; \{x_1,\ldots,x_n\})\leq C_0 \Big\}
		\,=\, O\big( n^{-2/3} \big).
	\end{displaymath}
}

\subsection{A partially linear model}
\label{S4+}

In this section we add a linear trend and a seasonal component to our original model (\ref{eq2.1}).
To simplify notation, we suppose again that the time points $x_1,\ldots,x_n$ are equidistant and that the
seasonal component has period~$p$.
This leads to the partially linear model 
\begin{equation}
	\label{eq5.1}
	Y_t \,=\, \big((x_1-\bar{x}_n)/(x_n-x_1) \big)\gamma_0^0 \,+\, \gamma_{(t \bmod p)+1}^0
	\,+\, m_0(x_t) \,+\, \varepsilon_t, \qquad t=1,\ldots,n.
\end{equation}
Throughout this section we assume that {\bf (A1)} and {\bf (A2)} are fulfilled.

In a first step, the parameter $\gamma^0=(\gamma_0^0,\gamma_1^0,\ldots,\gamma_p^0)^T$ is estimated by least squares.
We rewrite (\ref{eq5.1}) in vector/matrix form,
\begin{displaymath}
	Y \,=\, X\gamma^0 \,+\, \bar{m}_0 \,+\, \varepsilon,
\end{displaymath}
where $Y=(Y_1,\ldots,Y_n)^T$, $\bar{m}_0=(m_0(x_1),\ldots,m_0(x_n))^T$, $\varepsilon=(\varepsilon_1,\ldots,\varepsilon_n)^T$, and
\begin{displaymath}
	X \,=\, \left( \begin{array}{cccccc}
		\quad (x_1-\bar{x}_n)/(x_n-x_1) \quad & \quad 1 \quad & \quad 0 \quad & \quad \hdots \quad & \quad \hdots \quad & \quad 0 \quad \\
		\quad (x_2-\bar{x}_n)/(x_n-x_1) \quad & \quad 0 \quad & \quad 1 \quad & \quad 0 \quad & \quad \hdots \quad & \quad 0 \quad \\
		\quad \vdots \quad & \quad \vdots \quad & \quad \ddots \quad & \quad \ddots \quad & \quad \ddots \quad & \quad \vdots \quad \\
		\quad (x_{p-1}-\bar{x}_n)/(x_n-x_1) \quad & \quad 0 \quad & \quad \hdots \quad & \quad 0 \quad & \quad 1 \quad & \quad 0 \quad \\
		\quad (x_p-\bar{x}_n)/(x_n-x_1) \quad & \quad 0 \quad & \quad \hdots \quad & \quad \hdots \quad & \quad 0 \quad & \quad 1 \quad \\
		\quad (x_{p+1}-\bar{x}_n)/(x_n-x_1) \quad & \quad 1 \quad & \quad 0 \quad & \quad \hdots \quad & \quad \hdots \quad & \quad 0 \quad \\
		\quad (x_{p+2}-\bar{x}_n)/(x_n-x_1) \quad & \quad 0 \quad & \quad 1 \quad & \quad 0 \quad & \quad \hdots \quad & \quad 0 \quad \\
		\quad \vdots \quad & \quad \vdots \quad & \quad \ddots \quad & \quad \ddots \quad & \quad \ddots \quad & \quad \vdots \quad \\
		\quad \vdots \quad & \quad \vdots \quad & \quad \vdots \quad & \quad \vdots \quad & \quad \vdots \quad & \quad \vdots \quad
	\end{array} \right) ,
\end{displaymath}
with $\bar{x}_n=(x_1+\cdots +x_n)/n$.
To ensure identifiability of the parameters we choose $m_0$ such that 
$\sum_{t=1}^n t\, m_0(x_t)=0$ and
$\sum_{1\leq t\leq n\colon \, t \bmod p =k} m_0(x_t)=0$
for $k=1,\ldots,p$, which implies that
\begin{equation}
	\label{eq5.2}
	X^T \bar{m}_0 \,=\, 0_{p+1}.
\end{equation}
It is easy to see that
\begin{displaymath}
	n^{-1} X^TX \ninfty \mbox{Diag}\big[ 1/12, 1/p, \ldots, 1/p \big],
\end{displaymath}
which means that $X^TX$ is regular for sufficiently large~$n$ and
\begin{eqnarray*}
	\widehat{\gamma}_n & = & \argmin_\gamma \big\| Y \,-\, X\gamma \big\|^2 \\
	& = & \argmin_\gamma \big\| X\gamma^0 \,+\, \bar{m}_0 \,+\, \varepsilon \,-\, X\gamma \big\|^2 \\
	& = & \argmin_\gamma \big\| X\gamma^0 \,+\, \varepsilon \,-\, X\gamma \big\|^2 \\
	& = & \big( X^TX \big)^{-1} X^T\big( X\gamma^0 \,+\, \varepsilon \big).
\end{eqnarray*}
Then
\begin{eqnarray*}
	E \widehat{\gamma}_n & = & \gamma_0, \\
	E\big[ (\widehat{\gamma}_n-\gamma_0)(\widehat{\gamma}_n-\gamma_0)^T \big]
	& = & \big( X^TX \big)^{-1} X^T\Cov(\varepsilon)X \big( X^TX \big)^{-1} \,=\, O\big( n^{-1} \big).
\end{eqnarray*}

In a second step we estimate $m_0$ nonparametrically by wavelet thresholding.
Let
\begin{displaymath}
	\widetilde{Y}_t \,:=\, Y_t \,-\, \big( X\widehat{\gamma}_n \big)_t
	\,=\, m_0(x_t) \,+\, \big( X(\gamma^0-\widehat{\gamma}_n) \big)_t \,+\, \varepsilon_t.
\end{displaymath}
It follows from our identifiability condition (\ref{eq5.2}) that $\sum_{t=1}^n m_0(x_t)=0$.
Hence, $m_0$ can be represented as a linear combination of the wavelets and $\phi_0$ is not needed.
Let, for $(j,k)\in {\mathcal I}_n$,
\begin{eqnarray*}
	\widetilde{\beta}_{j,k} & = & \frac{1}{n} \sum_{t=1}^n \psi_{j,k}(x_t) \widetilde{Y}_t \\
	& = & \beta_{j,k}^0 \,+\, \frac{1}{n} \sum_{t=1}^n \psi_{j,k}(x_t) \big( X(\gamma^0-\widehat{\gamma}_n) \big)_t
	\,+\, \frac{1}{n} \sum_{t=1}^n \psi_{j,k}(x_t) \varepsilon_t.
\end{eqnarray*}
It follows from Lemma~\ref{LA.2} that $E\Big[ \big\| \widehat{\gamma}_n \,-\, \gamma^0 \big\|^4 \Big] \,=\, O\big( n^{-2} \big)$,
which implies
\begin{displaymath}
	\sup_t \,E\Big[ \big( X(\widehat{\gamma}_n-\gamma^0)_t \big)^4 \Big] \,=\, O\big( n^{-2} \big),
\end{displaymath}
and therefore, in conjunction with (\ref{eq4.11}),
\begin{displaymath}
	E\big[ \big( \widetilde{\beta}_{j,k} \,-\, \beta_{j,k}^0 \big)^4 \big] \,=\, O\big( n^{-2} \big) \qquad \forall (j,k)\in {\mathcal I}_n.
\end{displaymath}
This is an analog to equation (\ref{eq4.11}) which was the starting point for our calculations in the previous section
and we obtain the following result.

{\prop
	Suppose that  $(Y_t)_{t}$ satisfies \eqref{eq5.1} and that {\bf (A1)} and {\bf (A2)} are fulfilled. Then
	\begin{eqnarray*}
		E\big[ \|\widehat{\gamma}-\gamma_0\|^2 \big] & = & O\big( n^{-1} \big), \\
		\frac{1}{n} \sum_{t=1}^n \big( \widehat{m}_n(x_t) \,-\, m_0(x_t) \big)^2 & = & O_P\big( n^{-2/3}\, \log(n) \big).
	\end{eqnarray*}
}

\section{Simulations and data examples}
\label{S5}

\subsection{Simulations}\label{S51}

We illustrate the finite sample performance of the wavelet estimator with soft thresholding proposed in Section~\ref{SS4.2} using the following two trend functions:
$$
\begin{aligned}
	f(t)&=
	\begin{cases}
		1.5+t,\quad &t\in[0,\,1/2),\\
		0.1, \quad& t\in[1/2,\,2/3),\\
		3\sqrt{t-2/3}+0.1, \quad& t\in [2/3,1)
	\end{cases}\\
	g(t)&=
	\begin{cases}
		10t-\lfloor10t\rfloor,\quad &t\in[0,\,0.7),\\
		0.5, \quad& t\in[0.7,\,1).
	\end{cases}
\end{aligned}
$$

\begin{figure}[h]\label{f.functions}
	\includegraphics[width=10cm]{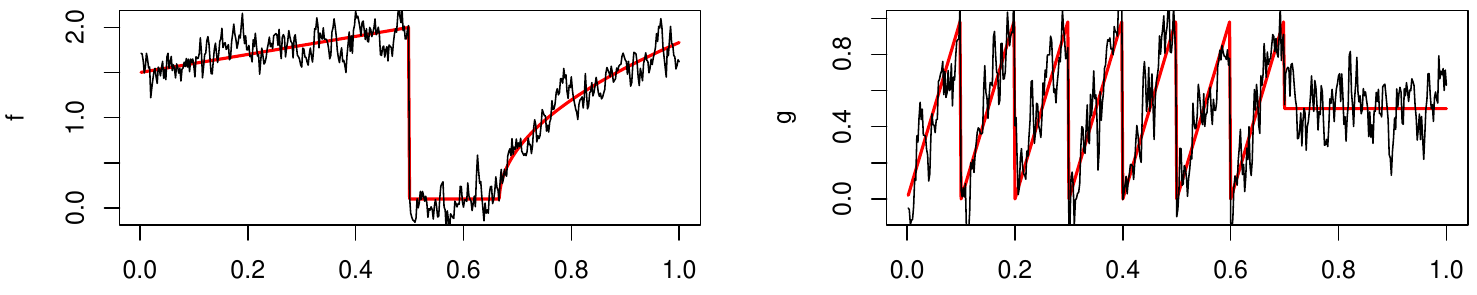}
	\caption{\textcolor{red}{\textbf{Red:}} function $f$ (left) and $g$ (right), \textbf{black:} data generated according to $Y_t=f(t/n)+\varepsilon_t$ (left) and $Y_t=g(t/n)+\varepsilon_t$ (right).}
\end{figure}

The shape of the first function is motivated by our data example displayed in Figure~\ref{fig.intro} and analyzed in Section~\ref{S52}. 
The second function is used to illustrate that our approach can successfully estimate rough functions with several jumps.
In both cases, we simulate an AR(1) noise process $(\varepsilon_t)_t$  with autoregressive parameter $a=0.7$ and i.i.d.~normal innovations with variance 0.01, see Fig.~\ref{f.functions}. 

\begin{figure}[h]
	\includegraphics[width=10cm]{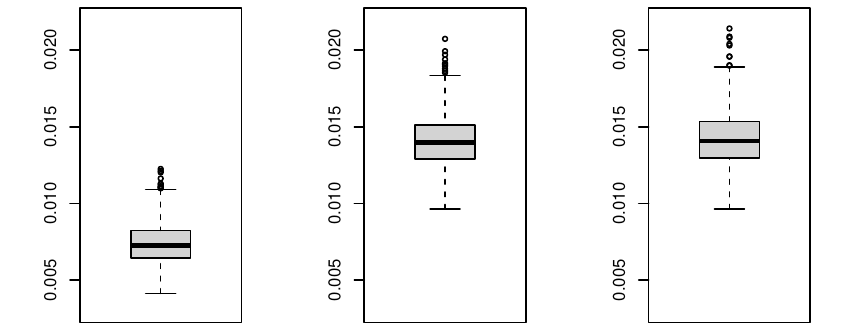}
	\caption{left: MSE of the wavelet estimator of $f$ with $K=0.1$, middle: MSE of the NW estimator of $f$ (rectangular kernel with $b=0.009$) right: MSE of the NW estimator of $f$ (Epanechnikov kernel with $b=0.007$). }\label{f.f-boxplot}
\end{figure}

\begin{figure}[h]
	\includegraphics[width=10cm]{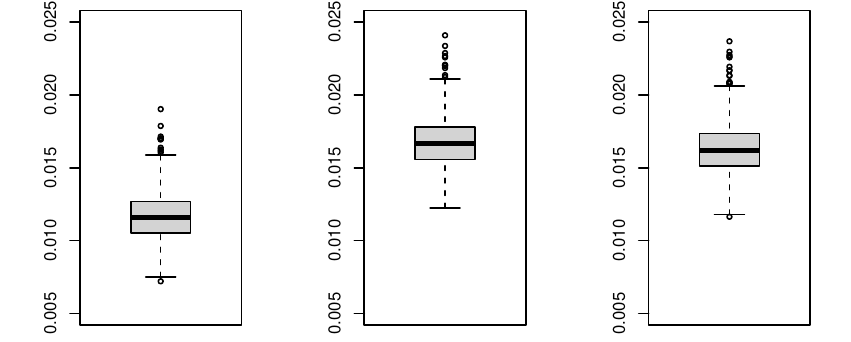}
	\caption{left: MSE of the wavelet estimator of $g$ with  $K=0.045$, middle: MSE of the Nadaraya-Watson estimator of $g$ (rectangular kernel with $b=0.006$) right: MSE of the Nadaraya-Watson estimator of $g$ (Epanechnikov kernel with  $b=0.007$).}\label{f.g-boxplot}
\end{figure}

We calculate the wavelet estimator proposed in Section~\ref{SS4.2} and compare it to the Nadaraya-Watson estimator generated with the function \texttt{kreg} from the \texttt{R} package \texttt{gplm}. For the Nadaraya-Watson estimator we consider both, the rectangular and the Epanechnikov kernel.  For our estimator as well as the kernel estimator we choose the tuning parameters (threshold parameter $K$ and bandwidth $b$, respectively) in an MSE-optimal manner using a grid search. Figures~\ref{f.f-boxplot} and \ref{f.g-boxplot} display the resulting box plots based on 1000 Monte Carlo iterations. 
In both settings, the  wavelet estimator outperforms the competing kernel estimators. Moreover, note that the MSE-optimal bandwidths for the kernel estimators are much smaller than the corresponding default values chosen by Scott's rule of thumb  (\citet[p.~152, eq.~6.42]{Sc92}) which are $b=0.186$ for the rectangular kernel and $b=0.145$ for the Epanechnikov kernel, respectively, in case of the function $g$. The latter bandwidths lead to oversmooth estimators of the trend function, while the choices in Figures~\ref{f.f-boxplot} and \ref{f.g-boxplot} result in an overfitting.

\subsection{A real data example}\label{S52}

\begin{figure}[H]
	\includegraphics[width=0.9\textwidth]{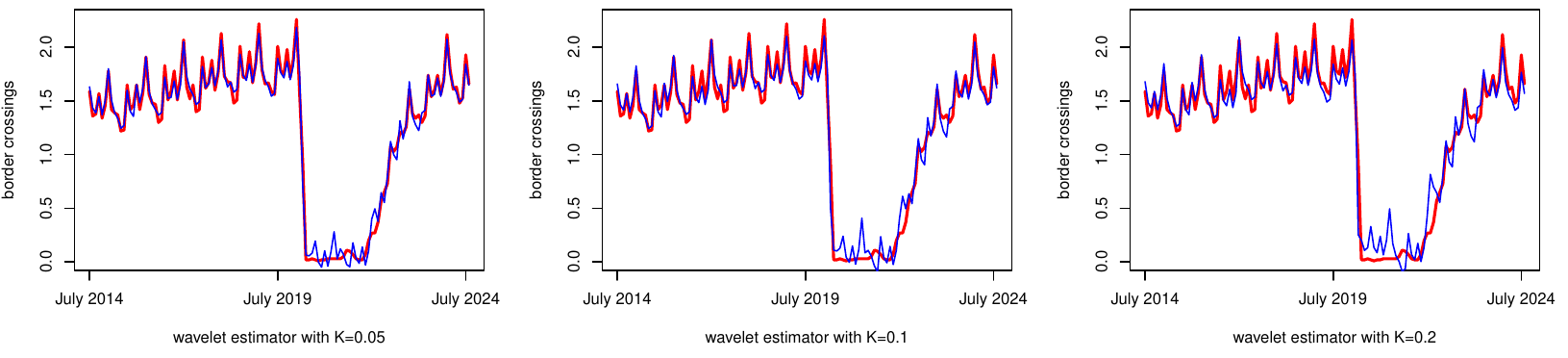}
	\caption{\textcolor{red}{red}: original data, \textcolor{blue}{blue}: (partially linear)  wavelet approximations with different choices of $K$ }\label{fig1}
\end{figure}

The data set contains monthly overseas arrival data in Australia. More precisely, it consists of  monthly recordings of international border crossings (in millions) from July 2014 to August 2024, retrieved from the Australian Bureau of Statistics\footnote{\url{https://www.abs.gov.au/statistics/industry/tourism-and-transport/overseas-arrivals-and-departures-australia/latest-release}}. As it can be seen from Figure~\ref{fig1}, there was an abrupt decay of the arrivals in April 2020 resulting from travel restrictions due to the COVID pandemic. 

We applied the partially linear model with soft thresholding, introduced in Section~\ref{S4+} with three different choices for the thresholding parameter $K=0.05,~K=0.1,~K=0.2$. For comparison, we additionally fitted a partially linear model, where the nonlinear part is estimated via classical kernel regression  (Nadaraya-Watson).  
From Figure~\ref{fig1} one can see that the wavelet-based estimator is perfectly able to adapt to the sharp decrease of the arrivals in 2020 irrespective of the choice of $K$.

\begin{figure}[H]
	\includegraphics[width=0.9\textwidth]{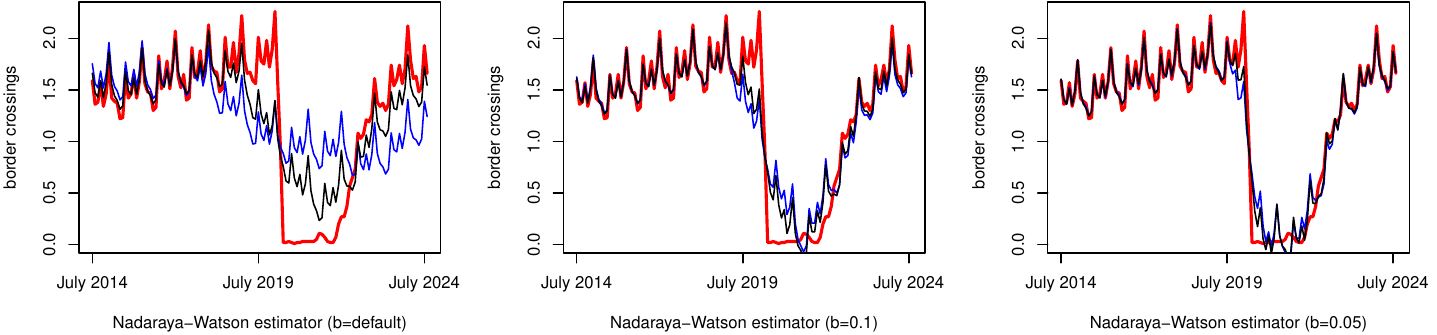}
	\caption{\textcolor{red}{red}: original data, \textcolor{blue}{blue}: (partially linear)  Nadaraya-Watson approximation with the rectangular kernel, \textcolor{black}{black}: (partially linear)  Nadaraya-Watson approximation with the Epanechnikov kernel. }\label{fig.data-nw}
\end{figure}

In contrast, the decay of the kernel-type estimator is clearly slower if the bandwidth is chosen via Scott's rule of thumb,
see~Figure~\ref{fig.data-nw} (left).  Reducing the bandwidth, we observe  that the estimator can adapt better to the rough path of the data. The question of choosing an appropriate bandwidth is much more important here than choosing the truncation parameter~$K$ for our wavelet estimator. The latter perfectly captures the COVID jump for all three choices of~$K$.

\section{Proofs of the main results}
\label{S6}

\begin{proof}[Proof of Proposition~\ref{P3.1}]
	Since the pairs $(\theta_1,\varepsilon_1),\ldots,(\theta_N,\varepsilon_N)$ are by assumption independent
	and identically distributed,
	it suffices to consider the Bayes risk for a single coefficient.
	Let $P_\theta$ be the distribution of~$Y_k$ given $\theta_k=\theta$ and let $T(Y_k)$ be an arbitrary
	estimator of~$\theta_k$. Then its Bayes risk w.r.t.~prior~$\pi$ and squared error loss is given by
	\begin{eqnarray*}
		\lefteqn{ E\big[ \big( T(Y_k) \,-\, \theta_k \big)^2 \big] } \\
		& = & \pi\big( \{0\} \big) \, 
		\Big\{ P_0\big( Y_k=-\lambda \big) T(-\lambda)^2 \,+\, P_0\big( Y_k=0 \big) T(0)^2 \,+\, P_0\big( Y_k=\lambda \big) T(\lambda)^2 \Big\} \\
		& & {} \,+\, \pi\big( \{-\lambda\} \big) \, 
		\Big\{ P_{-\lambda}\big( Y_k=-2\lambda \big) \big( T(-2\lambda)+\lambda \big)^2 
	\,+\, P_{-\lambda}\big( Y_k=-\lambda \big) \big( T(-\lambda)+\lambda \big)^2\\
			&&\hspace{2.5cm}\,+\, P_{-\lambda}\big( Y_k=0 \big) \big( T(0)+\lambda \big)^2 \Big\} \\
		& & {} \,+\, \pi\big( \{\lambda\} \big) \, 
		\Big\{ P_{\lambda}\big( Y_k=0 \big) \big( T(0)-\lambda \big)^2 \,+\, P_{\lambda}\big( Y_k=\lambda \big) \big( T(\lambda)-\lambda \big)^2\\
			&&\hspace{2.5cm}\,+\, P_{\lambda}\big( Y_k=2\lambda \big) \big( T(2\lambda)-\lambda \big)^2 \Big\} \\
		& = & (1-p) \, \big\{ (p/2)\, T(-\lambda)^2 \,+\, (1-p)\, T(0)^2 \,+\, (p/2)\, T(\lambda)^2 \big\} \\
		& & {} \,+\, (p/2)\, \big\{ (p/2)\, \big( T(-2\lambda)+\lambda \big)^2 \,+\, (1-p)\, \big( T(-\lambda)+\lambda \big)^2
		\,+\, (p/2)\, \big( T(0)+\lambda \big)^2 \big\} \\
		& & {} \,+\, (p/2)\, \big\{ (p/2)\, \big( T(0)-\lambda \big)^2 \,+\, (1-p)\, \big( T(\lambda)-\lambda \big)^2
		\,+\, (p/2)\, \big( T(2\lambda)-\lambda \big)^2 \big\} \\
		& = & (p/2)^2\; \big( T(-2\lambda)+\lambda \big)^2 \\
		& & {} \,+\, (1-p)\, (p/2)\; \big\{ T(-\lambda)^2 \,+\, \big( T(-\lambda)+\lambda \big)^2 \big\} \\
		& & {} \,+\, (1-p)^2\; T(0)^2 \,+\, (p/2)^2\; \big\{ \big( T(0)-\lambda \big)^2 \,+\, \big( T(0)+\lambda \big)^2 \big\} \\
		& & {} \,+\, (1-p)\, (p/2)\; \big\{ T(\lambda)^2 \,+\, \big( T(\lambda)-\lambda \big)^2 \big\} \\
		& & {} \,+\, (p/2)^2\; \big( T(2\lambda)-\lambda \big)^2.
	\end{eqnarray*}
	Therefore, the unique Bayes estimator $T^*(Y_k)$ is given by
	$T^*(-2\lambda)\,=\,-\lambda$, $T^*(-\lambda)\,=\,-\lambda/2$, $T^*(0)\,=\,0$, $T^*(\lambda)\,=\,\lambda/2$, $T^*(2\lambda)\,=\,\lambda$,
	and its Bayes risk is equal to
	\begin{eqnarray*}
		\lefteqn{ E\big[ \big( T^*(Y_k) \,-\, \theta_k \big)^2 \big] } \\
		& = & (p/2)^2\, 0 \\
		& & {} \,+\, (1-p) \, \big\{ (p/2)\, (\lambda/2)^2 \,+\, (p/2)\, (\lambda/2)^2 \big\} \\
		& & {} \,+\, (1-p)^2\, 0 \,+\, (p/2)\, \big\{ (p/2)\, \lambda^2 \,+\, (p/2)\, \lambda^2 \big\} \\
		& & {} \,+\, (1-p)\, \big\{ (p/2)\, (\lambda/2)^2 \,+\, (p/2)\, (\lambda/2)^2 \big\} \\
		& & {} \,+\, (p/2)^2\, 0 \\
		& = & 4\, (1-p)\, (p/2)\, (\lambda/2)^2 \,+\, 2\, (p/2)^2\, \lambda^2 \\
		& = & p\lambda^2/2 \,=\, \epsilon^2 q^{2/3}/2.
	\end{eqnarray*}
\end{proof}
\bigskip

\begin{proof}[Proof of Theorem~\ref{T3.1}]
	\begin{itemize}
		\item[(i)\quad]
		To prove (i) we distinguish between two cases.
		\begin{itemize}
			\item[1)\quad] $|\theta_k|< t/2$\\
			If $|\varepsilon_k|=|Y_k-\theta_k|\leq t/2$, then $|Y_k|< t$ and so $\widehat{\theta}_k=0$.
			Otherwise, we use the estimate
			$|\widehat{\theta}_k-\theta_k|\leq |\widehat{\theta}_k-Y_k|+|Y_k-\theta_k|\leq t+|\varepsilon_k|\leq 3|\varepsilon_k|$.
			Therefore,
			\begin{displaymath}
				E\big[ (\widehat{\theta}_k - \theta_k)^2 \big] \,\leq\, \theta_k^2 \,+\, 9\, E\big[ \varepsilon_k^2\, \1(|\varepsilon_k|>t/2) \big].
			\end{displaymath}
			\item[2)\quad] $|\theta_k|\geq t/2$\\
			In this case we use that $|\widehat{\theta}_k-\theta_k|\leq t\,+\, |\varepsilon_k|$, which implies by (\ref{eq3.11})
			\begin{displaymath}
				E\big[ (\widehat{\theta}_k - \theta_k)^2 \big] \,\leq\, 2\, t^2 \,+\, 2\,E\big[ \varepsilon_k^2 \big]
				\,\leq\, 2\, t^2 \,+\, 4\, \epsilon^2.
			\end{displaymath}
		\end{itemize}
		\item[(ii)\quad]
		We have
		\begin{eqnarray*}
			(1/N) \sum_{k\colon |\theta_k|\leq t/2} \min\big\{ \theta_k^2, t^2 \big\}
			& \leq & (t/2) \; (1/N) \sum_{k=1}^N |\theta_k| \\
			& = & (K\epsilon q^{-1/3}/2) \; \epsilon q \,=\, O\big( \epsilon^2 q^{2/3} \big).
		\end{eqnarray*}
		Moreover, it follows from Lemma~\ref{LA.1} that
		\begin{eqnarray*}
			E\big[ \varepsilon_k^2 \1(|\varepsilon_k|>t/2) \big]
			& = & 2\, \int_{t/2}^\infty \big( 1-Q_\varepsilon([-x,x]) \big)\, x \, dx \;+\; \big( 1-Q_\varepsilon([-t/2,t/2]) \big)\, (t/2)^2 \\
			& \leq & 2\, \int_{t/2}^\infty (\epsilon/x)^4 \, x\, dx \;+\; (\epsilon/(t/2))^4\, (t/2)^2 \,=\, 8\epsilon^4/t^2.
		\end{eqnarray*}
		This implies
		\begin{subequations}
			\begin{equation}
				\label{pt31.2a}
				\frac{1}{N}\sum_{k\colon |\theta_k|\leq t/2} E_\theta\big[ (\widehat{\theta}_k-\theta_k)^2 \big] \,=\, O\big( \epsilon^2 q^{2/3} \big).
			\end{equation}
			Since $\#\big\{k\in\{1,\ldots,N\}\colon \, |\theta_k|\geq t/2 \big\} \,\leq\, (1+\delta)\epsilon q/(t/2)$, we obtain
			\begin{eqnarray}
				\label{pt31.2b}
				\frac{1}{N}\sum_{k\colon |\theta_k|> t/2} E_\theta\big[ (\widehat{\theta}_k-\theta_k)^2 \big]
				& \leq & \big( 2\, t^2 + 4\, \epsilon^2 \big) \, \#\big\{k\in\{1,\ldots,N\}\colon \, |\theta_k|\geq t/2 \big\} \nonumber \\
				& \leq & \big( 2\, t^2 + 4\, \epsilon^2 \big) \, \epsilon q(1+\delta)/(t/2) \,=\, O\big( \epsilon^2 q^{2/3} \big).
			\end{eqnarray}
		\end{subequations}
		The second result follows from (\ref{pt31.2a}) and (\ref{pt31.2b}).
	\end{itemize}
\end{proof}

\begin{proof}[Proof of Theorem~\ref{T4.2}]
	We have by assumption that
	\begin{displaymath}
		E\big[ \big( \widetilde{\alpha}_0 \,-\, \alpha_0^0 \big)^2 \big]
		\,+\, \sum_{j=0}^{J_n^*-1} \sum_{k\colon\! (j,k)\in {\mathcal I}_n} E\big[ \big( \widetilde{\beta}_{j,k} \,-\, \beta_{j,k}^0 \big)^2 \big]
		\,=\, O\big( n^{-2/3} \big).
	\end{displaymath}
	Since $\var(\widetilde{\beta}_{j,k})=O(t_{n,j})$ for all $(j,k)\in{\mathcal I}_n$, $j\geq J_n^*$, we obtain
	from (i) of Theorem~\ref{T3.1} that
	\begin{eqnarray*}
		\lefteqn{ \sum_{j=J_n^*}^{J_n} \sum_{k\colon\! (j,k)\in {\mathcal I}_n} E\big[ (\widehat{\beta}_{j,k} - \beta_{j,k}^0)^2 \big] } \\
		& = & \sum_{j=J_n^*}^{J_n} \sum_{k\colon\! (j,k)\in {\mathcal I}_n} O\Big( \min\big\{ (\beta_{j,k}^0)^2, t_{n,j}^2 \big\}
		\;+\; \, E\big[ (\widetilde{\beta}_{j,k}-\beta_{j,k}^0)^2 \, \1\big( |\widetilde{\beta}_{j,k}-\beta_{j,k}^0|>t_{n,j}/2 \big) \big] \Big).
	\end{eqnarray*}
	Lemma~\ref{L4.1} reveals that
	\begin{displaymath}
		\sum_{j=J_n^*}^{J_n} \sum_{k\colon\! (j,k)\in {\mathcal I}_n}
		\min\big\{ (\beta_{j,k}^0)^2, t_{n,j}^2 \big\} \,=\, O\big( n^{-2/3} \big).
	\end{displaymath} 
	Let now $(j,k)\in {\mathcal I}_n$ with $j\geq J_n^*$. Then, since
	$P\big( \big| \widetilde{\beta}_{j,k}-\beta_{j,k}^0 \big| > x \big) \,\leq\, 1 \,-\, G(x/\sqrt{n})$ $\forall x\geq 0$,
	\begin{eqnarray*}
		\lefteqn{ E\Big[ \big( \widetilde{\beta}_{j,k} \,-\, \beta_{j,k}^0 \big)^2 \, \1\big( |\widetilde{\beta}_{j,k}-\beta_{j,k}^0| \,>\, t_{n,j}/2 \big)
			\Big] } \\
		& \leq & \int \big(n^{-1/2} x\big)^2 \, \1\big( |n^{-1/2}x| > K n^{-2/3} 2^{j/2-1} \big) \, dG(x) \\
		& = & n^{-1} \; \sum_{l=0}^\infty \int_{(Kn^{-1/6}2^{(j+l)/2-1}, Kn^{-1/6}2^{(j+l+1)/2-1}]} x^2 \, dG(x) \\
		& \leq & n^{-1} \; \sum_{l=0}^\infty \int_{(Kn^{-1/6}2^{(j+l)/2-1}, Kn^{-1/6}2^{(j+l+1)/2-1}]} \frac{x^4}{K^2n^{-1/3}2^{j+l-2}} \, dG(x) \\
		& = & 2^{-j}\, n^{-2/3}/K^2 \; \sum_{l=0}^\infty \int_{(Kn^{-1/6}2^{(j+l)/2-1}, Kn^{-1/6}2^{(j+l+1)/2-1}]} x^4/2^{l-2} \, dG(x).
	\end{eqnarray*}
	This implies
	\begin{eqnarray*}
		\lefteqn{ \sum_{j=J_n^*}^{J_n} \sum_{k\colon\! (j,k)\in {\mathcal I}_n} E\Big[ \big( \widetilde{\beta}_{j,k} \,-\, \beta_{j,k}^0 \big)^2 \,
			\1\big( |\widetilde{\beta}_{j,k}-\beta_{j,k}^0| \,>\, t_{n,j}/2 \big) \Big] } \\
		& \leq & \frac{4\, n^{-2/3}}{K^2} \sum_{m=0}^\infty \big( \frac{1}{2^0} \,+\, \cdots \,+\, \frac{1}{2^m} \big)
		\int_{(Kn^{-1/6}2^{m/2},Kn^{-1/6}2^{(m+1)/2}]} x^4\, dG(x) \\
		& \leq & \frac{8\, n^{-2/3}}{K^2} \int_{(0,\infty)} x^4 \, dG(x),
	\end{eqnarray*} 
	which completes the proof.                                                            
\end{proof}

\section{A few auxiliary results}
\label{S7}

{\lem 
	\label{LA.1}
	Let $F$ be the cumulative distribution function of a nonnegative random variable.
	Then, for any $t\geq 0$,
	\begin{displaymath}
		\int_t^\infty x^2 \, dF(x)
		\,=\, 2\, \int_t^\infty \big( 1-F(x) \big) x \, dx \;+\; \big( 1-F(t) \big) \, t^2.
	\end{displaymath}
}

\begin{proof}
	Let $b>t$ and
	\begin{displaymath}
		G(x) \,=\, \left\{ \begin{array}{ll} 0, & \quad \mbox{ if } x\leq 0, \\
			x^2, & \quad \mbox{ if } 0\leq x\leq b, \\
			b^2, & \quad \mbox{ if } x\geq b
		\end{array} \right. .
	\end{displaymath}
	Since $G$ is a continuous distribution function which has derivative $2x$ on $(0,b)$, we
	obtain by integration by parts
	\begin{eqnarray*}
		\int_t^b x^2 \, dF(x)
		& = & F(b)\, b^2 \;-\; F(t)\, t^2 \;-\; \int_t^b F(x)\, dG(x) \\
		& = & 2\, F(b) \int_0^b x\, dx \;-\; F(t)\, t^2 \;-\; 2\, \int_t^b F(x)\, x \, dx \\
		& = & 2\, \int_t^b \big( F(b)\,-\,F(x) \big) x \, dx \;+\; \big( F(b)\,-\,F(t) \big)\, t^2.
	\end{eqnarray*}
	The result follows with $b\to\infty$.
\end{proof}

{\lem
	\label{LA.2}
	Suppose that (A2) is fulfilled. Then
	\begin{displaymath}
		E\Big[ \big( \sum_{s=1}^n a_s\, \varepsilon_s \big)^4 \Big]
		\,\leq\, 3\, C_1^2\, \big( \sum_{s=1}^n a_s^2 \big)^2 \,+\, C_2\, \sum_{s=1}^n a_s^4.
	\end{displaymath}
}

\begin{proof}
	We have that
	\begin{eqnarray*}
		\lefteqn{ E\Big[ \big( \sum_{s=1}^n a_s\, \varepsilon_s \big)^4 \Big] } \\
		& = & \sum_{s,t,u,v=1}^n a_s\, a_t\, a_u\, a_v\; E[\varepsilon_s\, \varepsilon_t\, \varepsilon_u\, \varepsilon_v] \\
		& = & \sum_{s,t,u,v=1}^n a_s\, a_t\, a_u\, a_v\; \big\{ E[\varepsilon_s\, \varepsilon_t]\, E[\varepsilon_u\, \varepsilon_v]
		\,+\, E[\varepsilon_s\, \varepsilon_u]\, E[\varepsilon_t\, \varepsilon_v]
		\,+\, E[\varepsilon_s\, \varepsilon_v]\, E[\varepsilon_t\, \varepsilon_u] \big\} \\
		& & {} \,+\, \sum_{s,t,u,v=1}^n a_s\, a_t\, a_u\, a_v\; \cum(\varepsilon_s, \varepsilon_t, \varepsilon_u, \varepsilon_v) \\
		& =: & T_{n,1} \,+\, T_{n,2}.
	\end{eqnarray*} 
	Then 
	\begin{eqnarray*}
		T_{n,1} & = & 3\, \Big( \sum_{s,t=1}^n a_s\, a_t\; E[\varepsilon_s\, \varepsilon_t] \Big)^2 \\
		& \leq & 3\, \Big( \sum_{s,t=1}^n (a_s^2 + a_t^2)/2\; |\cov(\varepsilon_s, \varepsilon_t)| \Big)^2 \\
		& \leq & 3\; C_1^2\, \Big( \sum_{s=1}^n s_s^2 \Big)^2
	\end{eqnarray*}
	and since $|a_sa_ta_ua_v|\leq(a_s^4+a_t^4+a_u^4+a_v^4)/4$,
	\begin{eqnarray*}
		T_{n,2} & = & \sum_{s=1}^n a_s^4 \; \sum_{t,u,v=1}^n |\cum(\varepsilon_s, \varepsilon_t, \varepsilon_u, \varepsilon_v)| \,\leq \, C_2\, \sum_{s=1}^n a_s^4.
	\end{eqnarray*}
\end{proof}

\begin{proof}[Proof of Lemma~\ref{L4.1}]
	To simplify notation, let $K=1$, which means that $t_{n,j}=n^{-2/3}2^{j/2}$.
	Before we delve into details of the proof we consider an extreme case in order to provide some intuition.
	If $m_0=\1_{[c,x_n]}$ for any $c$ between $x_1$ and $x_n$, then we obtain that 
	at each scale $j$ only one of the coefficients~$\beta_{j,k}^0$ can be non-zero with $\beta_{j,k}^0=O(2^{-j/2})$.
	Since $2^{-j/2}\geq t_{n,j}$ if and only if $2^j\leq n^{2/3}$, we obtain that
	\begin{eqnarray*}
		S_n\big( m_0 \big)
		& \leq & \sum_{j\colon\! n^{1/3}\leq 2^j\leq n^{2/3}} t_{n,j}^2
		\,+\, \sum_{j\colon\! 2^j> n^{2/3}} \; \sum_{k\colon\! (j,k)\in {\mathcal I}_n} (\beta_{j,k}^0)^2 \\
		& = & O\bigg( n^{-4/3} \sum_{j\colon\! n^{1/3}\leq 2^j\leq n^{2/3}} 2^j \bigg) \,+\, O\bigg( \sum_{j\colon\! 2^j>n^{2/3}} 2^{-j} \bigg)
		\,=\, O\big( n^{-2/3} \big).
	\end{eqnarray*}
	In such a case, and in other cases as well, the desired result for $S_n(m_0)$ will be obtained by cutting the terms $|\beta_{j,k}^0|$
	at the corresponding thresholds $t_{n,j}$ up to a certain scale.
	
	Now we turn to the general case. We have that
	\begin{eqnarray}
		\label{pla3.1}                                                                                                               
		\big| \beta_{j,k}^0 \big| & \leq & (1/n)\, \sum_{t=1}^n \big| \psi_{j,k}^{(n)}(x_t) \big| \,
		\Big( \max\big\{ m_0(x_t)\colon x_t\in I_{j,k} \big\} \,-\, \min\big\{ m_0(x_t)\colon x_t\in I_{j,k} \big\} \Big) /2 \nonumber \\
		& \leq & 2^{-j/2} \, \mbox{TV}\big( m_0; I_{j,k} \big) \,=:\, c_{j,k}.
	\end{eqnarray}
	In order to estimate $S_n(m_0)$ we replace the terms $|\beta_{j,k}^0|$ by their upper estimates $c_{j,k}$.
	In what follows we make use of the
	relations between the $c_{j,k}$ at adjacent scales and of the nested structure of the intervals $I_{j,k}$.
	If $c_{j,k}>t_{n,j}$, which is equivalent to $\mbox{TV}(m_0; I_{j,k})>n^{-2/3}2^j$, then we obtain for $(j',k')\subseteq {\mathcal I}_n$
	with $I_{j,k}\subseteq I_{j',k'}$ that
	\begin{subequations}
		\begin{equation}
			\label{pla3.2a}
			c_{j',k'} \,=\, 2^{-j'/2} \mbox{TV}(m_0; I_{j',k'}) \,\geq\, 2^{-j'/2} \mbox{TV}(m_0; I_{j,k}) \,>\, n^{-2/3} 2^{j'/2} \,=\, t_{n,j'}.
		\end{equation}
		On the other hand, if $c_{j,k}\leq t_{n,j}$, which is equivalent to $\mbox{TV}(m_0; I_{j,k})\leq n^{-2/3}2^j$,
		then we obtain for $(j',k')\in {\mathcal I}_n$ with $I_{j',k'}\subseteq I_{j,k}$ that
		\begin{equation}
			\label{pla3.2b}
			c_{j',k'} \,=\, 2^{-j'/2} \mbox{TV}(m_0; I_{j',k'}) \,\leq\, 2^{-j'/2} \mbox{TV}(m_0; I_{j,k}) \,\leq\, n^{-2/3} 2^{j'/2} \,=\, t_{n,j'}.
		\end{equation}
	\end{subequations}
	Let ${\mathcal I}_{n,t}:=\big\{(j,k)\in {\mathcal I}_n\colon \; 2^j\geq n^{1/3}, \; c_{j,k}>t_{n,j}\big\}$ and
	${\mathcal I}_{n,c}:=\big\{(j,k)\in {\mathcal I}_n\colon\;  2^j\geq n^{1/3}, \; c_{j,k}\leq t_{n,j}\big\}$.
	Then
	\begin{equation}
		\label{pla3.3}
		S_n(m_0) \,\leq\, \sum_{(j,k)\in {\mathcal I}_{n,t}} t_{n,j}^2 \,+\, \sum_{(j,k)\in {\mathcal I}_{n,c}} c_{j,k}^2
		\,=:\, S_{n,t}(m_0) \,+\, S_{n,c}(m_0).
	\end{equation}
	In order to estimate the two terms on the right-hand side of (\ref{pla3.3}) we make use of the nested structure of the subsets
	${\mathcal I}_{n,t}$ and ${\mathcal I}_{n,c}$. 
	Let
	\begin{displaymath}
		{\mathcal I}_{n,t}^{top} \,:=\, \big\{(j,k)\in {\mathcal I}_{n,t}\colon\;
		(j',k')\not\in {\mathcal I}_{n,t} \mbox{ for all } I_{j',k'}\subseteq I_{j,k} \big\}
	\end{displaymath}
	be the collection of those intervals which are on top of the ``pyramid'' of intervals $I_{j,k}$ with $(j,k)\in {\mathcal I}_{n,t}$.
	Furthermore, let
	\begin{displaymath}
		{\mathcal I}_{n,c}^{bottom} \,:=\, \big\{(j,k)\in {\mathcal I}_{n,c}\colon\;
		I_{j-1,[k/2]}\in {\mathcal I}_{n,t} \big\}
	\end{displaymath}
	be the collection of those intervals from ${\mathcal I}_{n,c}$ which reside on top of an interval $I_{j,k}$ with $(j,k)\in {\mathcal I}_{n,t}$.
	
	We have, for $(j,k)\in {\mathcal I}_{n,t}^{top}$, that
	\begin{displaymath}
		t_{n,j}^2 \,<\, t_{n,j}\; 2^{-j/2}\mbox{TV}(m_0;I_{j,k}) \,=\, n^{-2/3}\, \mbox{TV}(m_0;I_{j,k}).
	\end{displaymath}
	Therefore, and since the intervals $I_{j,k}$ with $(j,k)\in {\mathcal I}_{n,t}^{(top)}$ are disjoint, we obtain that
	\begin{subequations}
		\begin{eqnarray}
			\label{pla3.4a}
			S_{n,t}(m_0)
			& = & \sum_{(j,k)\in {\mathcal I}_{n,t}^{top}} \; \sum_{(j',k')\in {\mathcal I}_n\colon\! I_{j,k}\subseteq I_{j',k'}} t_{n,j'}^2 \nonumber \\
			& \leq & \sum_{(j,k)\in {\mathcal I}_{n,t}^{top}} \; \sum_{j'\colon\! n^{1/3}\leq 2^{j'}\leq 2^j} t_{n,j'}^2 \nonumber \\
			& \leq & \sum_{(j,k)\in {\mathcal I}_{n,t}^{top}} 2\, n^{-2/3}\, \mbox{TV}(m_0;I_{j,k}) \nonumber \\
			& \leq & 2\, n^{-2/3}\, \mbox{TV}(m_0;(0,1]).
		\end{eqnarray}
		Note that we have for $(j,k)\in {\mathcal I}_{n,c}^{bottom}$ that $2^{j/2}c_{n,j}\leq n^{-2/3}$
		and for $(j',k')\in{\mathcal I}_n$ with $I_{j',k'}\subseteq I_{j,k}$ that $c_{j',k'}\leq 2^{(j-j')/2}c_{j,k}$.
		Therefore, and since the intervals $I_{j,k}$ with $(j,k)\in {\mathcal I}_{n,c}^{bottom}$ are disjoint, we obtain that
		\begin{eqnarray}
			\label{pla3.4b}
			S_{n,c}(m_0)
			& = & \sum_{(j,k)\in {\mathcal I}_{n,c}^{bottom}} \; \sum_{(j',k')\in {\mathcal I}_n\colon\!I_{j',k'}\subseteq I_{j,k}} c_{j',k'}^2 \nonumber \\
			& \leq & \sum_{(j,k)\in {\mathcal I}_{n,c}^{bottom}} \; \sum_{(j',k')\in {\mathcal I}_n\colon\! I_{j',k'}\subseteq I_{j,k}} 
			2^{(j-j')/2} c_{j,k} \; 2^{-j'/2}\mbox{TV}(m_0; I_{j',k'})  \nonumber \\
			& \leq & \sum_{(j,k)\in {\mathcal I}_{n,c}^{bottom}} 2^{j/2} c_{j,k} \sum_{j'\geq j} 2^{(j-j')} \mbox{TV}(f; I_{j,k}) \nonumber \\
			& \leq & 2\, n^{-2/3}\, \mbox{TV}(m_0; (0,1]).
		\end{eqnarray}
	\end{subequations}
	The result follows from (\ref{pla3.3}), (\ref{pla3.4a}), and (\ref{pla3.4b}).
\end{proof}

	%
	%
	
	\begin{funding}
		The second author was supported in part by the Oberfrankenstiftung (project: FP01054).
	\end{funding}

	\bibliographystyle{harvard}

\end{document}